\documentclass[review]{elsarticle}
\usepackage{amsmath,amsthm,amssymb,amsbsy,fancyhdr,graphicx,
            calc,ifthen,float,psfrag,mathrsfs,hhline}
\input epsf
\PassOptionsToPackage{ctagsplt,Righttag}{amstex}
\usepackage{mathtools}
\usepackage{framed}
\usepackage{dsfont}
\usepackage{ulem}
\usepackage{titlesec}
\newcounter{hours}
\newcounter{minutes}
\newcommand{\printtime}{\setcounter{hours}{\time/60}%
                        \setcounter{minutes}{\time-\value{hours}*60}%
\ifthenelse{\value{hours}<10}{0}{}\thehours:%
\ifthenelse{\value{minutes}<10}{0}{}\theminutes}
\usepackage{amsmath}
\usepackage{algorithm}
\usepackage{algorithmic}
\usepackage{amssymb}
\usepackage{bm}
\usepackage{mathtools}
\usepackage[font=scriptsize]{caption}
\usepackage[font=scriptsize]{subcaption}
\newcommand{\norm}[1]{\ensuremath{\left\| #1 \right\|}}

\newcommand{\pder}[2]{\ensuremath{\frac{\partial #1}{\partial #2}}} %1st partial derivative
\newcommand{\Ebb}{\ensuremath{\mathbb{E}}}

\newcommand{\Rbb}{\ensuremath{\mathbb{R} }}

\newcommand\Bbm{{\ensuremath{\bm{B}}}}

\newcommand\Fbm{{\ensuremath{\bm{F}}}}

\newcommand\Pbm{{\ensuremath{\bm{P}}}}

\newcommand\Xbm{{\ensuremath{\bm{X}}}}

\newcommand\cbm{{\ensuremath{\bm{c}}}}

\newcommand\fbm{{\ensuremath{\bm{f}}}}
\newcommand\gbm{{\ensuremath{\bm{g}}}}

\newcommand\nbm{{\ensuremath{\bm{n}}}}

\newcommand\tbm{{\ensuremath{\bm{t}}}}
\newcommand\ubm{{\ensuremath{\bm{u}}}}

\newcommand\xbm{{\ensuremath{\bm{x}}}}

\newcommand\lambdabold{{\ensuremath{\boldsymbol{\lambda}}}}

\newcommand\rhobold{{\ensuremath{\boldsymbol{\rho}}}}

\newcommand\xibold{{\ensuremath{\boldsymbol{\xi}}}}

\newcommand\Lambdabold{{\ensuremath{\boldsymbol{\Lambda}}}}

\newcommand\zerobold{\ensuremath{\mathbf{0}}}
\newcommand\onebold{\ensuremath{\mathbf{1}}}

\begin{document}
\bibliographystyle{plain}
\title{Design optimization in unilateral contact using pressure constraints and Bayesian optimization}
\author[1]{Jingyi Wang}
\author[1]{Jerome Solberg}
\author[1]{Mike A. Puso}
%\author[1,2]{Daniel A. Tortorelli}
\author[1]{Eric B. Chin}
\author[1]{Cosmin G. Petra}
\affiliation[1]{organization={Lawrence Livermore National Laboratory},country={USA}}
\affiliation[2]{organization={Department of Mechanical Sciences and Engineering, University of Illinois at Urbana-Champaign},
country={USA}}
%=============================================================================
%
\begin{abstract}
  Design optimization problems, \textit{e.g.}, shape optimization, that involve deformable bodies in unilateral contact are challenging as they require robust contact solvers, complex optimization methods that are typically gradient-based, and sensitivity derivations.
  Notably, the problems are nonsmooth, adding significant difficulty to the optimization process.
  We study design optimization problems in frictionless unilateral contact subject to pressure constraints, using both gradient-based and gradient-free optimization methods, namely Bayesian optimization.
  The contact simulation problem is solved via the mortar contact and finite element methods. 
  For the gradient-based method, we use the direct differentiation method to compute the sensitivities of the cost and constraint function with respect to the design variables.
  Then, we use Ipopt to solve the optimization problems.
  For the gradient-free approach, we use a constrained Bayesian optimization algorithm based on the standard Gaussian Process surrogate model. 
  We present numerical examples that control the contact pressure, inspired by real-life engineering applications, to demonstrate the effectiveness, strengths and shortcomings of both methods.
  Our results suggest that both optimization methods perform reasonably well for these nonsmooth problems. 
\end{abstract}

\maketitle
\section{Introduction}

Contacting bodies are ubiquitous in engineering problems, covering, \textit{e.g.}, crash worthiness test, part assembly, metal forming~\cite{laursen1992formulation}, etc.
Many of these problems are solved with finite element methods~\cite{wriggers2006computational}. 
However, combining design optimization and finite element methods  to \emph{design} contacting systems rather than \emph{analyze} them received little attention.
This is partially due to the complexity of combining state-of-the-art mortar contact solvers with optimization algorithms and the nonsmoothness that exists in such problems~\cite{hilding1999nonsmooth}.

Finite element contact solvers have evolved over the years and many techniques have been developed to solve contact problems in a robust and efficient manner~\cite{hallquist1985contact,puso2004contact}. Among them, the combined Uzawa-augmented Lagrangian method (ALM) has proven reliable. Recently, researchers have used the interior-point method, to solve frictionless contact problems~\cite{temizer2014}. 
The constraint function, \textit{i.e.}, gap function, can be computed using node-on-node, node-on-segment or mortar method.
The mortar segment-on-segment method is preferred because it passes the patch test and eliminates element locking~\cite{puso2004contact,puso2004contact2,zavarise2009node}.
In this work, we use the popular interior-point library Ipopt and the mortar method to solve the contact problem.

Previous shape and topology optimization work has applied gradient-based optimization methods to design systems with contact. 
In~\cite{fernandez2020topology}, the authors applied topology optimization to design nonlinear elastic bodies in contact using
Uzawa-ALM and mortar methods, adjoint sensitivities and Ipopt.
In~\cite{fancello2006}, the authors used the ALM method to deal with stress-based constraints and contact boundary conditions. The nonsmoothness of the solution to the contact problems is  highlighted as the reason for the design to be highly dependent on initial design and optimization algorithms.
In~\cite{petersson1997topology}, the design problems of sheets in contact are shown to be nonsmooth and solved approximately through a subgradient method.
In~\cite{bruggi2013stress}, a stress-based approach is presented for the design of truss-like elastic structures through imposing stress constraints.
Sensitivities and gradient-based methods are used in topology optimization problems in unilateral contact in~\cite{stromberg2010topology}.
In~\cite{pedersen2004crashworthiness}, the authors study crashworthiness optimization of frame structures using a topology optimization formulation and sensitivities.
Contact pressure distribution in plane stress with friction is optimized in~\cite{kristiansen2020} through the smoothing of a pressure function. While the study is limited to obstacle contact problems, the authors argued for the necessity of compliance constraint in addition to the objective function consisting of the mean and variance of pressure. 
In~\cite{niu2019topology}, the authors solved shape optimization with uniform contact pressure constraints. While the potential contact surface is not allowed to change, the actual contact region is affected by the shape design.
Though nonsmoothness caused by changing contact region is not discussed in the paper, the results show promise for  such problems in two dimension with method of moving asymptotes~\cite{svanberg1995mma}.
It is worth pointing out that most papers on design optimization in contact use node-on-segment contact or other simplified treatment of contact analysis, and not mortar contact formulation, rendering the methods less accurate. 
Moreover, design objective or constraints that are explicit in pressure are rare. 

Despite the success and advancement of these important works, there remain a few challenges that hinder the further development of a robust and efficient gradient-based optimizer.
The first one is the problem-dependent sensitivity computation based on complex contact solvers. Not only are derivation and implementation of sensitivities time-consuming, the necessary information to compute them is often difficult to obtain. For instance, many design variables change the shape and geometry of the bodies and thus the discretized mesh. The sensitivities of the mesh with respect to the design variables might be unavailable, particularly if third-party mesh tools are employed.
Second, the design optimization problem in contact can be nonsmooth, especially when considering explicit pressure objective and constraints~\cite{hilding1999nonsmooth}, leaving researchers with little convergence guarantee even with the most successful optimizers such as Ipopt or MMA~\cite{zillober1993mma}. 
Here, a smooth function means that it is continuously differentiable with Lipschitz first-order derivatives~\cite{Nocedal_book}.
During the optimization, the surface region in contact constantly changes with updated design variables, which could render the objective and constraint functions non-differentiable.   

To address the first challenge, we propose using the Bayesian optimization method for design problems in contact. 
Bayesian optimization is a gradient-free method that treats the objective and constraints as `black-box' functions that can be approximated with surrogate models, \textit{e.g.}, Gaussian Process (GP)~\cite{frazier2018}. 
Bayesian optimization has been widely applied to structural design~\cite{mathern2021}, process optimization of additive manufacturing~\cite{wang2023optimization}, inertial confinement fusion design~\cite{wang2024icf}, machine learning~\cite{brochu2010}, etc.
In~\cite{hsu2008insole}, a gradient-free subproblem approximation method is used for the optimization of an insole for lower plantar fascia stress.
To build the surrogate model, at each iteration, a sample point is collected by running the finite element simulation, where no sensitivity computation is necessary. 
Yet Bayesian optimization, particularly in the constrained case, has its own shortcomings. 
Without sensitivities, the algorithm is less likely to produce a local minimum that satisfies first-order optimality conditions~\cite{Nocedal_book}. In addition, for large scale problems, while the sampling can be easily parallelized, the surrogate model, \textit{e.g.}, Gaussian Processes (GPs) can become computationally expensive to train with increasing number of samples.

To study the second challenge,  we derive, design and implement sensitivities for problems with pressure objective and constraints, and solve them using the well-established interior-point solver Ipopt. 
The sensitivities are computed based on the differentiation of the Karush–Kuhn–Tucker (KKT) systems of the contact solutions.  
While no new algorithm is proposed to provide convergence theory, we demonstrate the validity of the gradient-based approach through numerical examples with changing contact region and pressure constraints that are inspired by real-life applications.
Further, we compare both gradient-based and gradient-free approaches on our numerical examples. Our main contributions are
\begin{itemize}
  \item derivation and implementation of a gradient-based method using sensitivities for design optimization in unilateral contact with pressure constraints
  \item Design and implementation of constrained Bayesian optimization on design optimization in unilateral contact
  \item Novel design and insights produced by our approach in numerical examples. 
\end{itemize}

The paper is organized as follows. In section~\ref{se:intro}, we briefly discuss the finite element contact mechanics models in the optimization form. In section~\ref{se:design}, we formulate the design optimization problem and its sensitivities based on the finite element model. We propose the constrained Bayesian optimization method in section~\ref{se:cbo}. The numerical examples are presented in section~\ref{se:exp}. Finally, the conclusion is given in section~\ref{se:conclusion}.

\section{Contact problem in optimization formulation}\label{se:intro}
By way of introduction, consider a solid body $B$ in the Euclidean space $\Ebb^n$ and let $\Xbm$ denote the location of its material points in the undeformed/reference configuration $\Omega$ with boundary $\partial \Omega$. 
The body is subjected to an external load field $f_{ext}$ and as such
the same material points are displaced to locations in the deformed configuration $\omega$ such that $\xbm=\Xbm+\ubm(\Xbm)$, where $\ubm:\Omega\to \Ebb^3$. The deformation gradient tensor $\Fbm $ is defined as $\Fbm=\onebold+\pder{\ubm}{\Xbm}$, where $\onebold$ is the mixed identity tensor.
We assume $B$ undergoes infinitesimal deformations and is comprised of an elastic material with internal energy function $\phi_{e}(\Fbm,\xbm)$. We also assume the external force $f_{ext}$ is conservative such that it is derivable from an external energy function $\phi_{ext}(\ubm,\Xbm)$.

Now consider two distinct bodies $\Omega^{(1)}$ and $\Omega^{(2)}$ with material points $\Xbm^{(i)}, i=1,2$ that are potentially in contact. The internal and external energies for $\Omega^{(i)}$ are $\phi^{(i)}_e(\Fbm^{(i)},\Xbm^{(i)})$ and $\phi_{ext}^{(i)}(\ubm^{(i)},\Xbm^{(i)})$.

\begin{figure}
  \centering
  \includegraphics[width=0.85\textwidth]{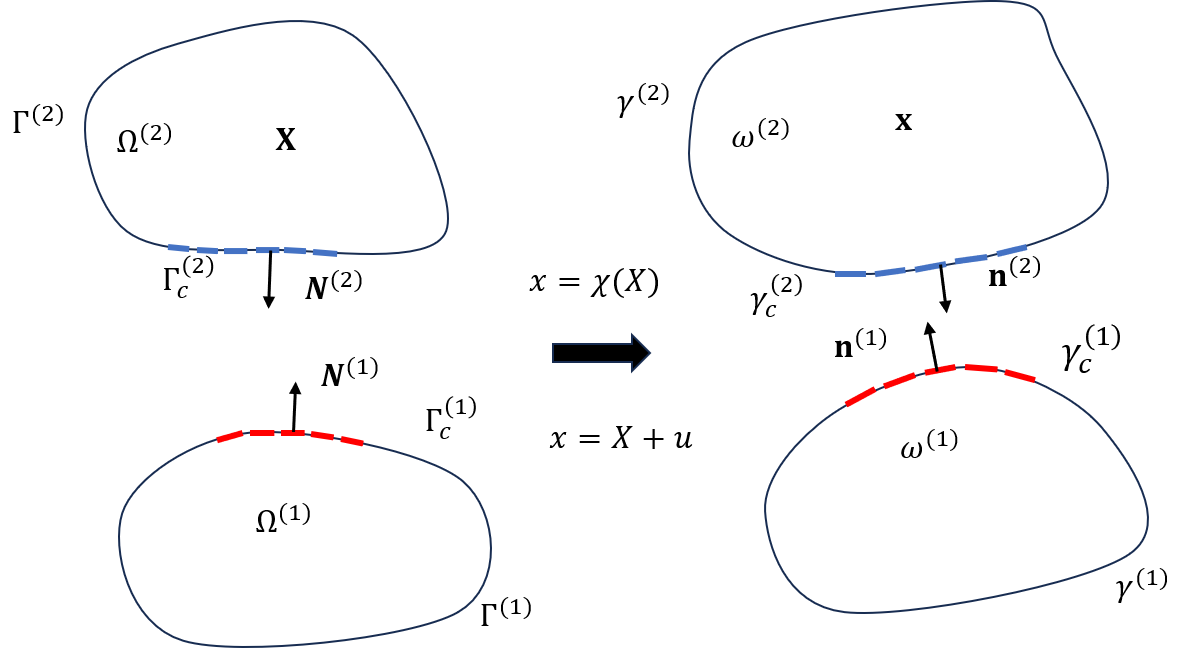}
	\caption{Reference and deformed configurations of bodies in potential contact.}
\label{fig:bodies}
\end{figure}

The boundary $\partial \Omega^{(i)}$ of $\Omega^{(i)}$ is decomposed into the complementary surfaces $\Gamma_u^i,\Gamma_t^i$ and $\Gamma_c^i$ of Figure~\ref{fig:bodies}.
Displacement is prescribed over $\Gamma_u^{(i)}$ and $\ubm^{(i)}(\Xbm^{(i)})= \bar{\ubm}^{(i)}(\Xbm^{(i)})$. 
The surface traction is prescribed over $\Gamma_t^{(i)}$: $\tbm^{(i)}(\Xbm^{(i)}) = \bar{\tbm}^{(i)}(\Xbm^{(i)})$. The surface $\Gamma_c^{(i)}$ denotes regions that will potentially contact with the other body. 
In the deformed configurations $\omega^{(i)}$ with boundary $\partial \omega^{(i)}$, these three surfaces occupy regions
    $\partial \omega^{(i)}$, $\gamma_u^{(i)}$, $\gamma_t^{(i)}$ and $\gamma_c^{(i)}$, respectively.

The unilateral contact constraint is based on the principle of impenetrability of matter, \textit{i.e.},
\begin{equation} \label{def:cur-boundary-c}
 \centering
  \begin{aligned}
     \omega^{(1)} \cap \omega^{(2)} = \emptyset.
  \end{aligned}
\end{equation}
We impose~\eqref{def:cur-boundary-c} via the gap function $g:\gamma_c^{(2)}\to \Rbb$ and constraint  
\begin{equation} \label{def:gap-function-node}
 \centering
  \begin{aligned}
     g(\xbm^{(2)}) \geq 0, \ \forall \xbm^{(2)} \in  \gamma_c^{(2)}.
  \end{aligned}
\end{equation}
To define $g$, we assume that each location $\xbm^{(2)}\in \gamma_c^{(2)}$ has a unique outward unit normal $\nbm^{(2)}(\xbm^{(2)})$.
The projection along $\nbm^{(2)}(\xbm^{(2)})$  is used to identify a potential contact point $\xbm^{(1)}\in\gamma_c^{(1)}$.
The signed distance between these points define the gap, \textit{i.e.}, $g=\nbm^{(2)}\cdot(\xbm^{(1)}-\xbm^{(2)})$~\cite{sewerin2020,puso2004contact}. 
The constraint~\eqref{def:gap-function-node} is one-sided and can bias the analysis~\cite{solberg2007}, but in general is effective and widely adopted.

To simplify the ensuing presentation, we define $\ubm = (\ubm^{(1)},\ubm^{(2)})$, $\Xbm = (\Xbm^{(1)},\Xbm^{(2)})$ and $\xbm = (\xbm^{(1)},\xbm^{(2)})$. 
As such, the frictionless contact for elastic bodies can be posed as an optimization problem:  
\begin{equation} \label{eqn:opt-contact}
 \centering
  \begin{aligned}
   &\underset{\substack{\ubm}}{\text{minimize}} 
	  & & f(\Fbm, \ubm,\Xbm) \\
   &\text{subject to}
	  & & g(\xbm) \geq 0, \ \forall \xbm \in \gamma_c^{(2)},\\
  \end{aligned}
\end{equation}
where $f$ is the total energy, \textit{i.e.},
\begin{equation} \label{eqn:opt-objective}
 \centering
  \begin{aligned}
     f(\Fbm,\ubm,\Xbm) = \sum_{i=1}^2 \int_{\Omega^{(i)}} \phi^{(i)}_e(\Fbm^{(i)},\Xbm^{(i)}) dv+\int_{\Omega^{(i)}}\phi_{ext}^{(i)}(\ubm^{(i)},\Xbm^{(i)}). 
  \end{aligned}
\end{equation}

We apply standard Galerkin finite element discretization to where in the displacement field is discretized with piecewise smooth interpolation functions $\Phi^{(i)}_j:\Omega^{(i)}\to \Omega$, local coordinates $\xibold\in\Rbb^3$, nodal displacement vectors $\ubm_{h,j}^{(i)}$ and $\ubm_{h,j}$
such that
\begin{equation} \label{eqn:discretization-1}
 \centering
  \begin{aligned}
     \ubm(\xibold,\ubm_{h,j}) = \sum_{j=1}^{n}  \Phi_j(\xibold) \ubm_{h,j}, 
  \end{aligned}
\end{equation}
where $n$ is the number of nodes. The subscript $h$ is used to denote nodal vectors and other discretized quantities, \textit{e.g.}, $\gamma_h^{(i)},i=1,2$. 
The discretized energy objective for an elastic material is  
\begin{equation} \label{eqn:opt-objective-dist}
 \centering
  \begin{aligned}
     f_h({\ubm}_h) = \sum_{i=1}^2  \frac{1}{2}\ubm_h^{(i),T} K^{(i)} \ubm_h^{(i)}  - \fbm_{ext}^{(i),T} \ubm_h^{(i)} 
         = \frac{1}{2} \ubm_h^{T} K \ubm_h  - \fbm_{ext}^T \ubm_h. 
  \end{aligned}
\end{equation}
where $K$ is the stiffness matrix and an integrated quantity, and $\fbm_{ext}$ is the external force vector. 

There are multiple ways to discretize and compute the gap function.
We choose the mortar method 
which integrates the pointwise gap function $g$ on $\gamma_{c,h}^{(2)}$.
Let $S$ be the index set of all nodes on the discretized contact boundary $\gamma_{c,h}^{(2)}$, the mortar gap function is defined as a vector function $\gbm_h:\Rbb^n\to\Rbb^m, m= |S|$, with components
\begin{equation} \label{eqn:opt-gap-dist}
 \centering
  \begin{aligned}
    \gbm_{h,j}(\ubm_h) = \int_{\gamma_{c,h}^{(2)}} \Phi_j(\xibold) g(\xibold,\ubm_h) d\xibold,
  \end{aligned}
\end{equation}
for all $j\in S$. Note that the local coordinates reduce to two dimensions on the surface.
Therefore, the discretized optimization problem of~\eqref{eqn:opt-contact} with mortar contact formulation is 
\begin{equation} \label{eqn:opt-contact-dist}
 \centering
  \begin{aligned}
   &\underset{\substack{\ubm_h\in\Rbb^n}}{\text{minimize}} 
	  & & f_h(\ubm_h) \\
   &\text{subject to}
	  & & \gbm_h(\ubm_h) \geq \zerobold.\\
  \end{aligned}
\end{equation}

The KKT conditions of~\eqref{eqn:opt-contact-dist} is that there exist Lagrange multipliers $\lambdabold_h \in\Rbb^m$ such that  
\begin{equation} \label{eqn:KKT-1}
 \centering
  \begin{aligned}
	  &\nabla f_h(\ubm_h) - \lambdabold_h \nabla \gbm_h(\ubm_h)  = 0\\
	  &\gbm_h(\ubm_h) \geq 0,  \ \lambdabold_h \geq 0, \ \lambda_{h,i} g_{h,i} = 0, \ i=1,...,m.
  \end{aligned}
\end{equation}
Equations~\eqref{eqn:KKT-1} comprise the standard contact finite element problem~\cite{eck2005unilateral} and 
is the basis for further deriving sensitivities needed for gradient-based design optimization.
The Lagrange multiplier $\lambdabold_h$ is often interpreted as the pressure on the contact surface~\cite{wriggers2006computational}. 
%The interior-point way of solving~\eqref{eqn:opt-contact} is through
%\begin{equation} \label{eqn:opt-contact-ip}
% \centering
%  \begin{aligned}
%   &\underset{\substack{u\in \Rbb^n}}{\text{minimize}} 
%	  & & f(u)\\
%   &\text{subject to}
%	  & & g(u) =s,\\
%       &&& s \geq 0,
%  \end{aligned}
%\end{equation}
%and then relax the complementarity condition: 
%\begin{equation} \label{eqn:KKT-1}
% \centering
%  \begin{aligned}
%	  &\nabla f(u) - \lambda \nabla g(u)  = 0\\
%	  &g(u) \geq 0,  \ \lambda \geq 0, \ \lambda_i g_i = \mu e, \ i=1,...,m.
%  \end{aligned}
%\end{equation}
%where $\mu \to 0$.
\section{Design optimization and sensitivities}\label{se:design}
In this section, we consider design optimization problems that are subjected pressure constraints.
Let $\rhobold\in \Rbb^p$ be the design variable. The mathematical form of our target problems is 
\begin{equation} \label{eqn:opt-1st}
 \centering
  \begin{aligned}
   &\underset{\substack{\rhobold}}{\text{minimize}} 
	  & & a(\rhobold,\bar{\ubm},\bar{\lambdabold})\\
   &\text{subject to}
	  & & \cbm(\rhobold,\bar{\ubm},\bar{\lambdabold}) \geq \zerobold,\\
  \end{aligned}
\end{equation}
where $a(\rhobold,\cdot,\cdot):\Rbb^p\times\Rbb^n\times\Rbb^m\to\Rbb$ and $\cbm(\rhobold,\cdot,\cdot):\Rbb^p\times\Rbb^n\times\Rbb^m\to\Rbb^q$ are the objective and constraint functions, respectively.
The displacement $\bar{\ubm}$ and contact pressure $\bar{\lambdabold}$ are obtained from the solutions to a contact mechanics problem~\eqref{eqn:opt-contact}. 

Applying the same finite element discretization,
the parameterized and discretized forms of~\eqref{eqn:opt-1st} can be written as 
\begin{equation} \label{eqn:opt-1st-dist}
 \centering
  \begin{aligned}
   &\underset{\substack{\rhobold}}{\text{minimize}} 
	  & & a_h(\rhobold,\bar{\ubm}_h,\bar{\lambdabold}_h)\\
   &\text{subject to}
	  & & \cbm_h(\rhobold,\bar{\ubm}_h,\bar{\lambdabold}_h) \geq \zerobold, \\
  \end{aligned}
\end{equation}
where given $\rhobold$, 
\begin{equation} \label{eqn:opt-2nd-dist}
 \centering
  \begin{aligned}
 \bar{\ubm}_h = \  &\underset{\substack{\ubm_h}}{\text{argmin}} \ 
	  & & f_h(\ubm_h,\rhobold) \\
   &\text{subject to}
	  & & \gbm_h(\ubm_h,\rhobold) \geq \zerobold,\\
  \end{aligned}
\end{equation}
and $\bar{\lambdabold}_h$ is the Lagrange multiplier associated with $\bar{\ubm}_h$. 
In many engineering applications such as shape optimization, $\Xbm_h$ can be a function of $\rhobold$, which results in~\eqref{eqn:opt-2nd-dist} being dependent on $\Xbm_h$ implicitly.

In order to use gradient-based optimization methods to solve~\eqref{eqn:opt-1st-dist}, we need to compute the state sensitivities against $\rhobold$ given by 
\begin{equation} \label{eqn:sens-1}
 \centering
  \begin{aligned}
      \pder{\bar{\ubm}_h}{\rhobold}, \ \pder{\bar{\lambdabold}_h}{\rhobold}.
  \end{aligned}
\end{equation}
For a low-dimensional design space where the dimension $p$ is small, a primal method suffices, while adjoint methods can be used when $p$ is large.
Denote the derivative of $f_h(\cdot,\cdot)$ and $\gbm_h(\cdot,\cdot)$ with respect to the $i$th variable by $\partial_i f_h(\cdot,\cdot)$ and $\partial_i \gbm_h(\cdot,\cdot)$, respectively, where $i=1,2$.
The solutions $\bar{\ubm}_h,\bar{\lambdabold}_h$ satisfy the parameterized KKT conditions
\begin{equation} \label{eqn:KKT-rhobold}
 \centering
  \begin{aligned}
	  &\partial_1 f_h(\bar{\ubm}_h,\rhobold) - \bar{\lambdabold}_h \partial_1 \gbm_h(\bar{\ubm}_h,\rhobold)  = \zerobold,\\
	  &\gbm_h(\bar{\ubm}_h,\rhobold) \geq \zerobold,  \ \bar{\lambdabold}_h \geq \zerobold, \ \bar{\Lambdabold}_{h} \gbm_{h}(\bar{\ubm}_h,\rhobold) = \zerobold, 
  \end{aligned}
\end{equation}
where $\bar{\Lambdabold}_h$ is a diagonal matrix whose values are the vector $\bar{\lambdabold}_h$.

Next, denoting the second-order derivatives of a function by $\partial^2_{i,j},i,j=1,2$ and differentiating the optimality conditions~\eqref{eqn:KKT-rhobold} with $\rhobold$, we have 
\begin{equation} \label{eqn:sens-2}
 \centering
  \begin{aligned}
	  \partial_{1,1}^2 f_h(\bar{\ubm}_h,\rhobold) \pder{\bar{\ubm}_h}{\rhobold} + \partial_{1,2}^2 f_h(\bar{\ubm}_h,\rhobold)
           -  [\partial_{1} \gbm_h(\bar{\ubm}_h,\rhobold)]^T \pder {\bar{\lambdabold}_h}{\rhobold} -& \\
           \pder{\bar{\ubm}_h}{\rhobold}^T [\partial_{1,1}^2 \gbm_h(\bar{\ubm}_h,\rhobold)]^T \bar{\lambdabold}_h 
           - [\partial_{1,2}^2 \gbm_h(\bar{\ubm}_h,\rhobold)]^T \bar{\lambdabold}_h= \zerobold&,\\
	   \pder{ \bar{\Lambdabold}_{h}}{\rhobold} \gbm_{h}(\bar{\ubm}_h,\rhobold) + \bar{\Lambdabold}_{h} \partial_1 \gbm_h(\bar{\ubm}_h,\rhobold) \pder{\bar{\ubm}_h}{\rhobold}+\bar{\Lambdabold}_{h} \partial_2 \gbm_{h}(\bar{\ubm}_h,\rhobold) = \zerobold&.
  \end{aligned}
\end{equation}
Solving the linear system of~\eqref{eqn:sens-2} gives us~\eqref{eqn:sens-1}.
We can then solve~\eqref{eqn:opt-1st} with gradient-based optimization methods by treating $\bar{\ubm}_h$ and $\bar{\lambdabold}_h$ as implicit functions of $\rhobold$. The gradients and Jacobian of the objective and constraint functions with respect to $\rhobold$ are computed with 
\begin{equation} \label{eqn:sens-3}
 \centering
  \begin{aligned}
    &\frac{d a_h}{d\rhobold} = \pder{a_h}{\rhobold} + \pder{a_h}{\bar{\ubm}_h}\pder{\bar{\ubm}_h}{\rhobold}+ \pder{a_h}{\bar{\lambdabold}_h}\pder{\bar{\lambdabold}_h}{\rhobold},\\
    &\frac{d \cbm_h}{d\rhobold} = \pder{\cbm_h}{\rhobold} + \pder{\cbm_h}{\bar{\ubm}_h}\pder{\bar{\ubm}_h}{\rhobold}+ \pder{\cbm_h}{\bar{\lambdabold}_h}\pder{\bar{\lambdabold}_h}{\rhobold}.
  \end{aligned}
\end{equation}
The second-order derivatives such as $\partial_{1,2}^2 f_h(\bar{\ubm}_h,\rhobold)$ are arrived at directly through the chain rule.
In many cases, $\rhobold$ could affect $\Xbm_h$ and the calculation of the second-order derivatives of the gap function is obtained again through chain rule. 

To solve the contact problem~\eqref{eqn:opt-2nd-dist}, we employ the interior-point solver Ipopt and HiOp~\cite{ipopt,petra2018hiop}, also referred to as the forward solver, with the full analytical Hessian. 
To solve the design optimization problem~\eqref{eqn:opt-1st-dist} using gradient-based optimization methods, the sensitivities are computed by solving the linear system~\eqref{eqn:sens-2} and the derivatives in~\eqref{eqn:sens-3} are calculated. 
For a robust and efficient outcome, we use Ipopt with its Hessian approximation option. The gradient-based optimization is implemented in MATLAB.

\section{Constrained Bayesian optimization}\label{se:cbo}
For the gradient-free optimization approach, we design and apply constrained Bayesian optimization with Gaussian Process (GP). 
Bayesian optimization is an iterative black-box optimization method that does not require the analytic forms of the objective functions or their derivatives. At each iteration, GP surrogate models replace both the objective and constraint functions and an acquisition function is used to select new sample points. Applied to design optimization, each sample point requires a solution to the parameterized contact problem~\eqref{eqn:opt-2nd-dist}. The optimization variables are typically bounded by a closed and bounded set $\rhobold\in C$.

A Gaussian process surrogate model assumes a multivariate joint Gaussian distribution between the variables $\rhobold$ and an output function $a_h(\rhobold)$, written as   
\begin{equation} \label{eqn:GP-1}
 \centering
  \begin{aligned}
  \begin{bmatrix}
   a_h(\rhobold_1) \\
   \vdots\\
           a_h(\rhobold_T)
  \end{bmatrix} \sim \mathcal{N}\left( 
      \begin{bmatrix}    
                 m(\rhobold_1)\\
 \vdots\\
 m(\rhobold_T)
      \end{bmatrix},
              \begin{bmatrix}
               k(\rhobold_1,\rhobold_1) \dots k(\rhobold_1,\rhobold_T)\\
      \vdots\\
               k(\rhobold_T,\rhobold_1) \dots k(\rhobold_T,\rhobold_T)
      \end{bmatrix}\ 
      \right). 
\end{aligned}
\end{equation}
Here, $\mathcal{N}$ is the normal distribution with $T$ samples. The design variable samples are denoted as $\rhobold_1,\ldots,\rhobold_T$. 
The function $m:\Rbb^p\to\Rbb$ is a user-chosen mean function for the distribution. The function $k:\Rbb^p\times\Rbb^p\to\Rbb$ denotes the covariance function. 
The posterior Gaussian probability distribution at a new sample point~$\rhobold$ can be inferred using Bayes' rule:
\begin{equation} \label{eqn:GP-post}
 \centering
  \begin{aligned}
  &a_h(\rhobold) | a_h(\rhobold_{1:T}) \sim \mathcal{N} (\mu(\rhobold),\sigma^2(\rhobold))\\
  &\mu(\rhobold)\ =\ k(\rhobold,\rhobold_{1:T}) k(\rhobold_{1:T},\rhobold_{1:T})^{-1} \left(a_h(\rhobold_{1:T})-m(\rhobold_{1:T}) \right) + m(\rhobold_{1:T})\\
  &\sigma^2(\rhobold)\ =\
k(\rhobold,\rhobold)-k(\rhobold,\rhobold_{1:T})k(\rhobold_{1:T},\rhobold_{1:T})^{-1}\sigma(\rhobold_{1:T},\rhobold)\ ,
\end{aligned}
\end{equation}
where the vector~$\rhobold_{1:T}$ is the notation for $[\rhobold_1,\dots,\rhobold_T]$ and 
\begin{equation} \label{eqn:GP-post-2}
 \centering
  \begin{aligned}
   k(\rhobold_{1:T},\rhobold_{1:T}) = [k(\rhobold_1),\dots,k(\rhobold_T); \dots; k(\rhobold_T,\rhobold_1),\dots,k(\rhobold_T,\rhobold_T)]. 
\end{aligned}
\end{equation}
The functions $\mu$ and $\sigma^2$ are the posterior mean and variance, respectively.
It is common to set the mean function $m$ to constant $0$~\cite{frazier2018}.
The covariance functions, or kernels, of the GP have significant impact on the accuracy of the surrogate model. 
A widely adopted kernel, the Squared Exponential Covariance Function, also  
known as the power exponential kernel, is defined as 
\begin{equation} \label{eqn:kernel-sec}
 \centering
  \begin{aligned}
  k(\rhobold,\rhobold';\theta)\ =\
\text{exp}\left(-\frac{\norm{\rhobold-\rhobold'}^2}{\theta^2}\right)\ ,
  \end{aligned}
\end{equation}
where~$\theta$ denotes the hyperparameters of the kernel. The hyperparameters of the GP surrogate model 
 often are optimized during the training or fitting step by maximizing the 
log-marginal-likelihood with an optimization algorithm. 

The choice of the acquisition function is critical to an efficient and effective Bayesian optimization algorithm.
If the acquisition function emphasizes on exploration, then typically sample points are chosen to reduce variance of the model. 
Otherwise, sample points are chosen to improve objective function value. 
One popular acquisition function is the expected improvement (EI)~\cite{brochu2010}, whose form can be written as
\begin{equation} \label{eqn:EI-1}
 \centering
  \begin{aligned}
       EI(\rhobold)\ =\ \begin{cases}
       (\mu(\rhobold) - a_h(\rhobold^+) - \xi)\Phi(Z)+\sigma(x)\phi(Z), &\text{if} \ \sigma(\rhobold)>0,\\
       0, &\text{if} \ \sigma(\rhobold) = 0,    
                 \end{cases}
  \end{aligned}
\end{equation}
where 
\begin{equation} \label{eqn:EI-2}
 \centering
  \begin{aligned}
       Z\ =\ \begin{cases}
       \frac{\mu(\rhobold)-a_h(\rhobold^+)-\xi}{\sigma(\rhobold)}, &\text{if} \ \sigma(\rhobold)>0,\\
       0, &\text{if} \ \sigma(\rhobold) = 0\ .    
                 \end{cases}
  \end{aligned}
\end{equation}
In the unconstrained setting, the variable~$\rhobold^+$ in~\eqref{eqn:EI-1} is defined as 
\begin{equation} \label{eqn:EI-3}
 \centering
  \begin{aligned}
	  \rhobold^+= \underset{\substack{ \rhobold_i\in \rhobold_{1:T}}}{\text{argmax}} \ a_h(\rhobold_i)
  \end{aligned}
\end{equation}
of existing $T$ samples, where $C$ is the bounds on $\rhobold$.
The trade-off parameter~$\xi\geq 0$ affects the trade-off between global search and local optimization~\cite{lizotte2008}. 
The functions $\phi$ and~$\Phi$ are the probability density function (PDF) and cumulative 
distribution function (CDF) of the standard normal distribution, respectively. 
The next sample point is chosen by maximizing the acquisition function. 
At the $k$th iteration, this leads to  
\begin{equation} \label{eqn:acquisition-1}
 \centering
  \begin{aligned}
	  \rhobold_{k} =  &\underset{\substack{\rhobold\in C}}{\text{argmax}} \ EI_k (\rhobold),
  \end{aligned}
\end{equation}
where $EI_k$ is the $EI$ function using samples accumulated up to the $k$th iteration.
In order to solve~\eqref{eqn:acquisition-1}, optimization algorithms including L-BFGS or random 
search can be used.

For constrained Bayesian optimization, the choices of acquisition functions are limited. 
Fortunately, $EI$ can be extended to constrained expected improvement $EI_C$, by incorporating feasibility penalties to the improvement function~\cite{gardner2014}.
The $\rhobold^+$ in~\eqref{eqn:EI-3} is defined as the feasible optimal solution from the $T$ samples 
\begin{equation} \label{eqn:EI-4}
 \centering
  \begin{aligned}
	  \rhobold^+=  &\underset{\substack{\rhobold_i\in \rhobold_{1:T}\\ \cbm_h(\rhobold_i)\geq 0 }}{\text{argmax}} \ a_h (\rhobold_i).
  \end{aligned}
\end{equation}
The expected constrained improvement acquisition function is 
\begin{equation} \label{eqn:EI-5}
 \centering
  \begin{aligned}
	  EI_C(\rhobold) = PF(\rhobold) EI(\rhobold),
  \end{aligned}
\end{equation}
where $PF(\rhobold) = Pr [c_h(\rhobold)\geq 0] $ is the probability of the constraints being satisfied.
For multiple inequality constraints, each constraint function is modeled as an independent GP surrogate model and is often simplified to be conditionally independent given $\rhobold$. 
That is, the probability of a feasible $\rhobold$ for constraints $\cbm_h(\rhobold) \geq 0$ satisfies 
\begin{equation} \label{eqn:EI-6}
 \centering
  \begin{aligned}
	  Pr(\cbm_{h}(\rhobold)\geq 0) = \prod_{i=1}^{q} Pr(\cbm_{h,i}(\rhobold)\geq 0).
  \end{aligned}
\end{equation}

The constrained Bayesian optimization algorithm is given in Algorithm~\ref{alg:bo}. 
\begin{algorithm}
 \caption{Constrained Bayesian optimization with GP}\label{alg:bo}
  \begin{algorithmic}[1]
	  \STATE{Choose initial sampling data $\rhobold_0$ in the bounds.} 
	  \STATE{Train GP surrogate models on the initial data.}
  \FOR{$i=0,1,2,\dots$}
	  \STATE{Evaluate acquisition function and find $\rhobold_i = \text{argmax}_{\rhobold} EI_{C,k}(x)$ within the bounds.}
	  \STATE{Run finite element simulation with design variables $\rhobold_i$.}
	  \STATE{Compute the objective $a_h(\rhobold_i)$ and constraints $\cbm_h(\rhobold_i)$ based on the simulation result. \;}
	  \STATE{Retrain the GP surrogate models with the addition of the new sample $\rhobold_i$, $a_h(\rhobold_i)$ and $\cbm_h(\rhobold_i)$.\;}
	  \STATE{Solve~\eqref{eqn:opt-1st-dist} as needed with the updated surrogate model.
	  Evaluate the stopping criteria. Exit if satisfied. \;}
  \ENDFOR
  \end{algorithmic}
\end{algorithm}
The training of the surrogate model refers to finding the optimal hyperparameters given the sample points. 
Line $8$ of Algorithm~\ref{alg:bo} can either return the feasible sample with the smallest objective or with the smallest posterior mean. 
Line $4$ can be solved by interior-point L-BFGS methods, our default choice, or random search algorithms.  
One stopping criteria is that the algorithm exists upon reaching a prescribed maximum number of iterations.

Algorithm~\ref{alg:bo} is implemented in MATLAB. We note that Bayesian optimization algorithms are often coded as maximization algorithms.  
In such case, $-a_h$ is the maximization objective.

\section{Numerical examples}\label{se:exp}
To demonstrate and compare the effectiveness of the gradient-based and gradient-free optimization methods, 
we present two numerical examples that are inspired by real-life engineering problems. 
These model problems can be easily extended to more complicated and larger-scale problems using the same optimization method we presented in this paper.
\subsection{Wedge problem}
\begin{figure}
  \centering
  \includegraphics[width=0.95\textwidth]{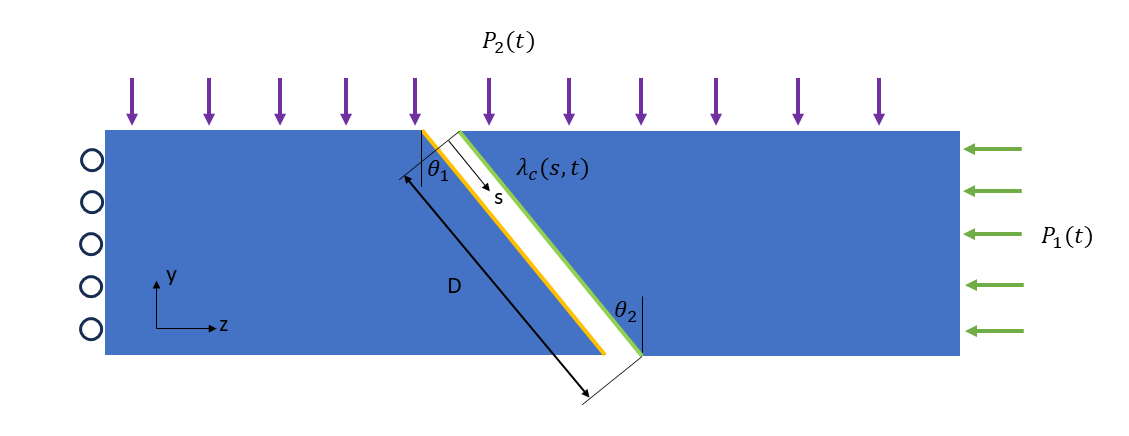}
	\caption{Wedge problem in the $y-z$ plane. $P_1(t)$ and $P_2(t)$ are two time-dependent loads.}
\label{fig:wedge2d}
\end{figure}

\begin{figure}
  \centering
  \includegraphics[width=0.75\textwidth]{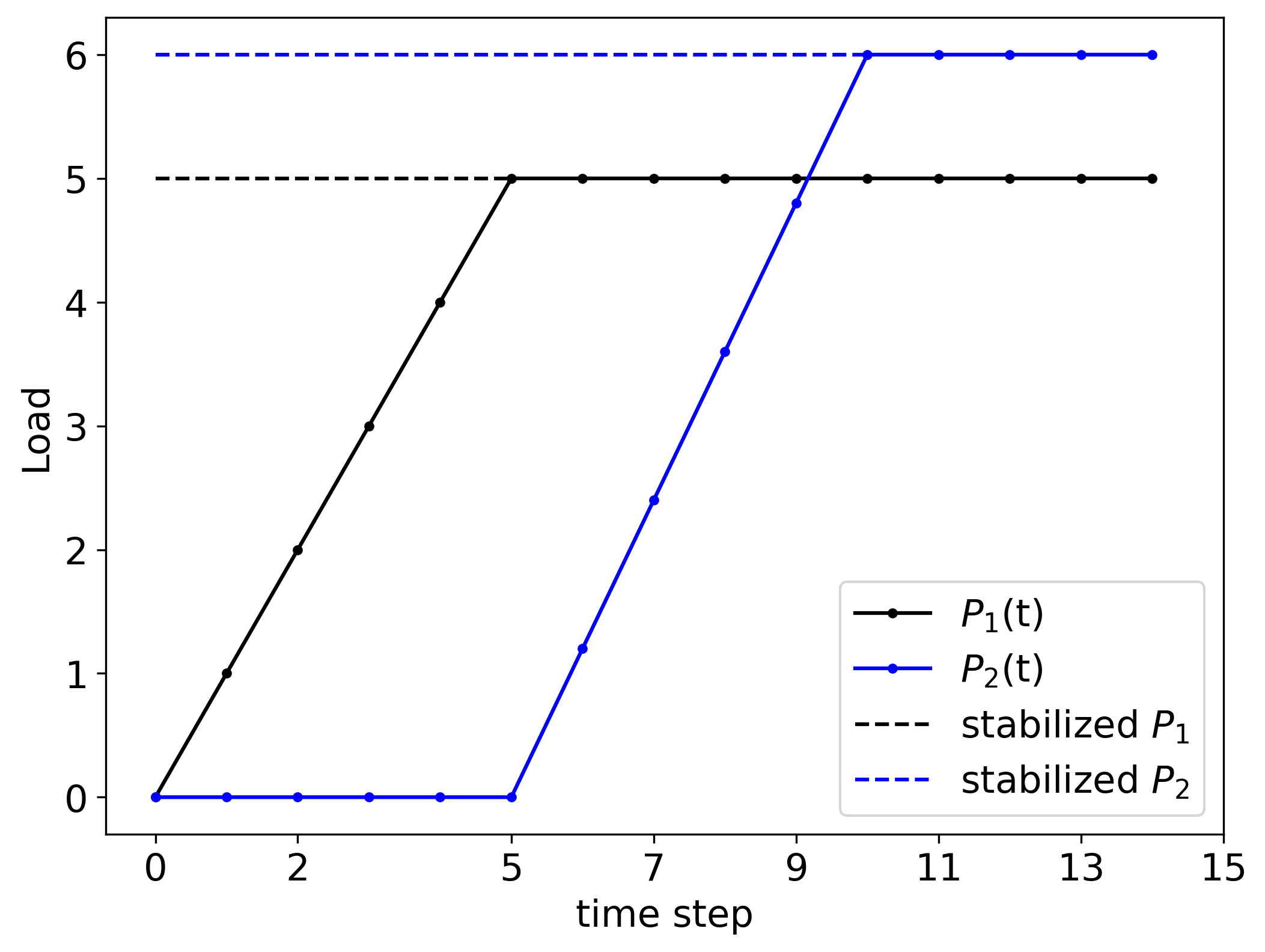}
	\caption{The pressure loads $P_1(t)$ and $P_2(t)$ as they change over time. The scalars $P_1$ and $P_2$ are the stabilized values of the loads.}
\label{fig:p1p2}
\end{figure}

The first numerical example is derived from the problem of designing high-current joints, whose illustration is given in Figure~\ref{fig:wedge2d}.  In this case, we consider a high-current joint where the current and magnetic field diffuse from the outside (positive z) inwards.  
At very early times, the current is concentrated at the very outer fiber of the joint, so it is important large contact pressure be maintained in the very beginning of the contact region, \textit{i.e.}, $\bar{\lambdabold}_h > \lambdabold_l$, so that arcing and/or joint separation dot occur.  As the current builds up, it diffuses into the surface to a depth $s=\alpha D$, but it is still important that the contact pressure be maintained above a minimum value $\lambdabold_l$ at the outer fiber and all the way down to the depth $s=\alpha D$, since every additional increment in current diffuses in from the outside.  At the same time, as current builds, a magnetic pressure is created, represented by $P_2(t)$ in Figure~\ref{fig:wedge2d}, which pushes inwards from the outside.  
The general design goal is to create a joint geometry which maintains adequate contact pressure within the joint, but without exceeding any material limits (e.g. locally $\bar{\lambdabold}_h < \lambdabold_m$).  At the optimum design, the amount of stabilized preload, represented by $P_1$ in Figure~\ref{fig:p1p2}, should be a minimum, as large values of contact preload require additional structure surrounding the joint.  
The joint is modeled in plane strain, but the extension to 3-D or axisymmetric analysis is trivial.
Based on the description above, the reference configuration in the $y-z$ plane can be seen in Figure~\ref{fig:wedge2d}.
The two external loads are depicted in Figure~\ref{fig:p1p2} below.

\begin{figure}
  \centering
  \includegraphics[width=\textwidth,trim={5cm 15cm 5cm 15cm},clip] {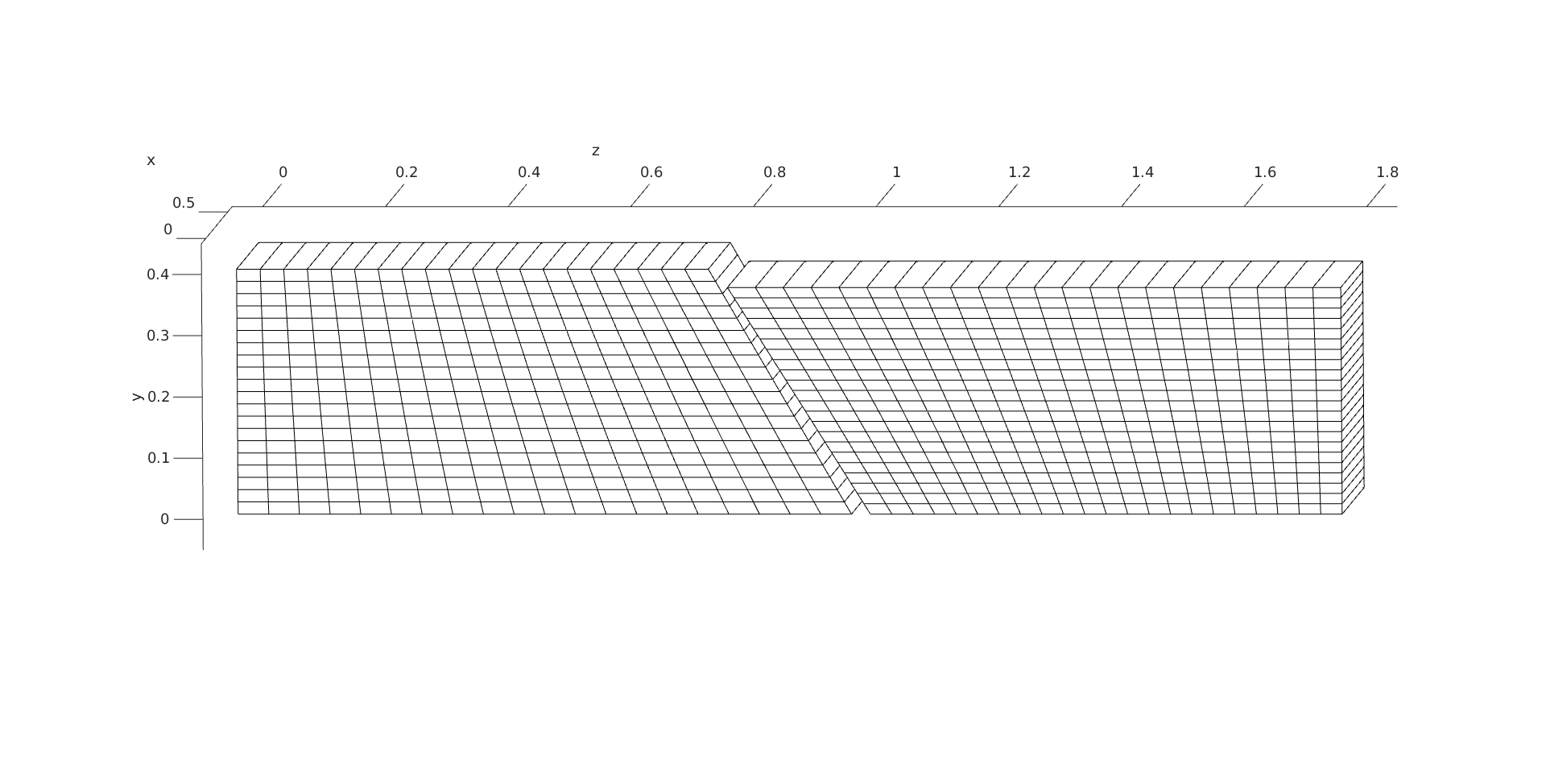}
	\caption{The reference configurations of the wedges in the initial design.} 
\label{fig:wedge0}
\end{figure}

The two wedges are subjected to two external loads in the form of pressure, where the first pressure $P_1(t)$ ramps up linearly from time $0$ and value $0$ to $5$ seconds and value $P_1$.
The second pressure $P_2(t)$ ramps up linearly from $5$ to $10$ seconds for a value from $0$ to $2$. 
The wedges are restricted in their movement in the $x$-direction and the left wedge has zero $z$-displacement on its left-most surface. 
The material for both wedges are chosen to be linear isotropic elastic material with the elastic modulus $E=200$ and Poisson's ratio $\nu=0.3$. For demonstration purpose, the problem is run unitless and can be easily extended to more complicated materials.
Both wedges are discretized with linear hexahedral elements, as shown in Figure~\ref{fig:wedge0}.
Given the plane strain boundary conditions imposed, the number of elements in the $x$ direction is reduced to $1$ to improve computational efficiency. The contact surface on the right wedge is chosen to be the surface for the mortar contact integrals.

As mentioned earlier, the design objective the wedge problem is the stabilized pressure $P_1$, while we impose constraints on the contact pressure on the upper part of the contact surfaces at both $t=5$ and $t=10$ seconds. 
The mathematical expression of the problem is
\begin{equation} \label{eqn:obj-1}
 \centering
  \begin{aligned}
   &\underset{\substack{\rhobold \in \Rbb^3}}{\text{minimize}} 
	  & & P_1 \\
   &\text{subject to}
	  && \lambda_{i}^e \geq \lambda_l , \ s\leq \alpha D,\\
	   &&& \lambda_{m} \geq \lambda_{i}^e , \ \forall s,
  \end{aligned}
\end{equation}
where the parameters are set to $\lambda_m=20$, $\lambda_l = 1.0$ and $\alpha = 0.4$ so that the top part of the contact surface are required to be in active contact. To reduce the potential nonsmoothness of pressure, the constraints are imposed on element segment pressure $\lambda_i^e$ for segment $i$, computed through the average of nodal pressure.
For the given discretization, we choose an element segment made of two elements, with a total of $11$ element segments, $22$ elements and $23$ nodes on the contact surface of the right wedge. Therefore, the constraints in~\eqref{eqn:obj-1} are imposed on a total of four element segments, equivalent to  eight elements, towards the top of the right wedge. 
The optimization variables are the two angles of the wedges $\theta_1$, $\theta_2$ and the pressure $P_1$. 
All of the three variables have bound constraints to ensure contact would occur at the two time steps.
Specifically, the bound constraints are 
\begin{equation} \label{eqn:ex1-bound}
 \centering
  \begin{aligned}
	   &30^{\circ} \leq \theta_1,\theta_2\leq 60^{\circ},\\
	   &0.5\leq P_1 \leq 1.5.
  \end{aligned}
\end{equation}

\begin{figure}
	\centering
	\begin{subfigure}{\linewidth}
        \centering
         \includegraphics[width=\textwidth,trim={5cm 15cm 5cm 15cm},clip]{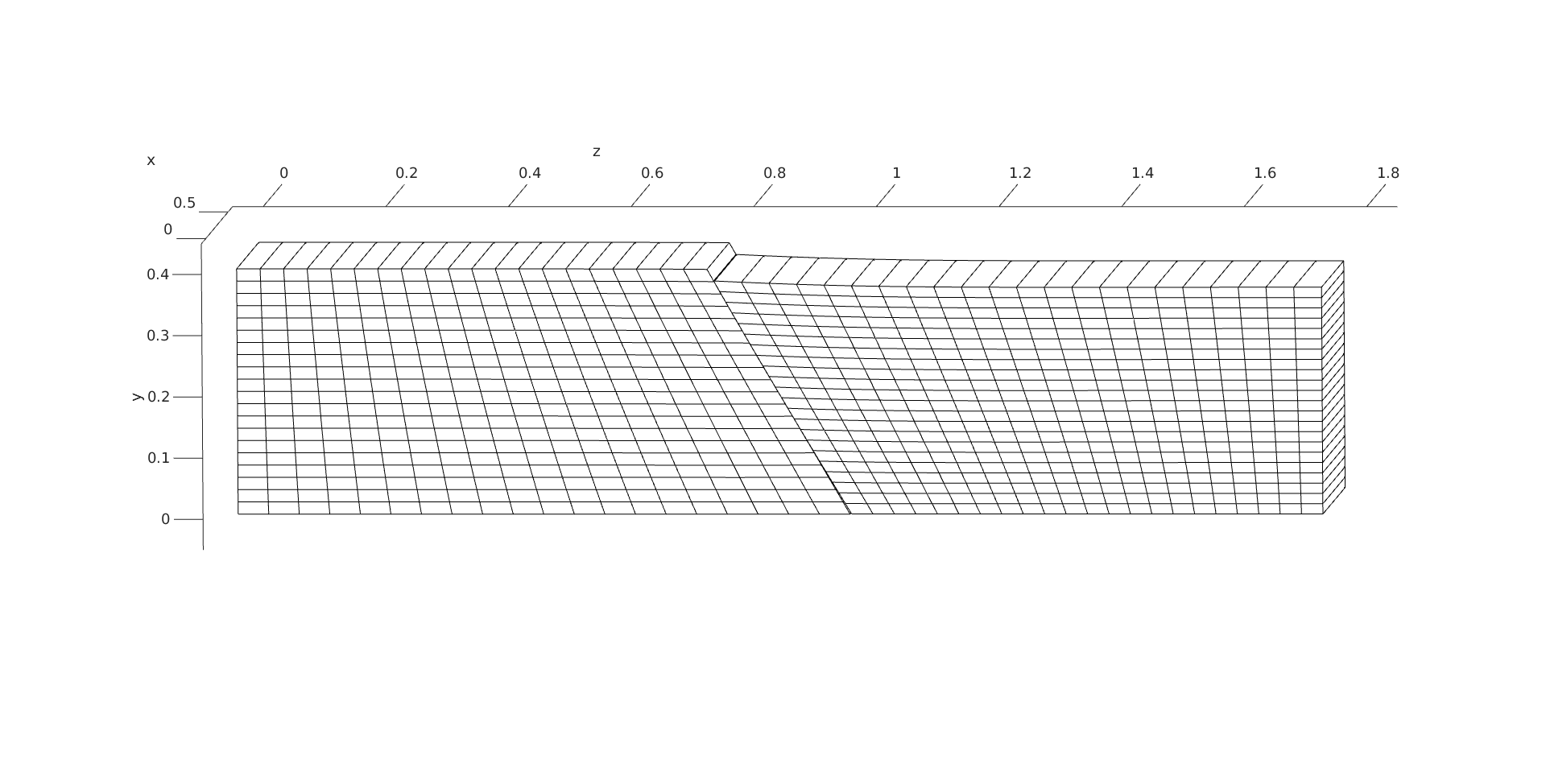}
     \end{subfigure} \\
	\begin{subfigure}{\linewidth}
        \centering
         \includegraphics[width=\textwidth,trim={5cm 15cm 5cm 15cm},clip]{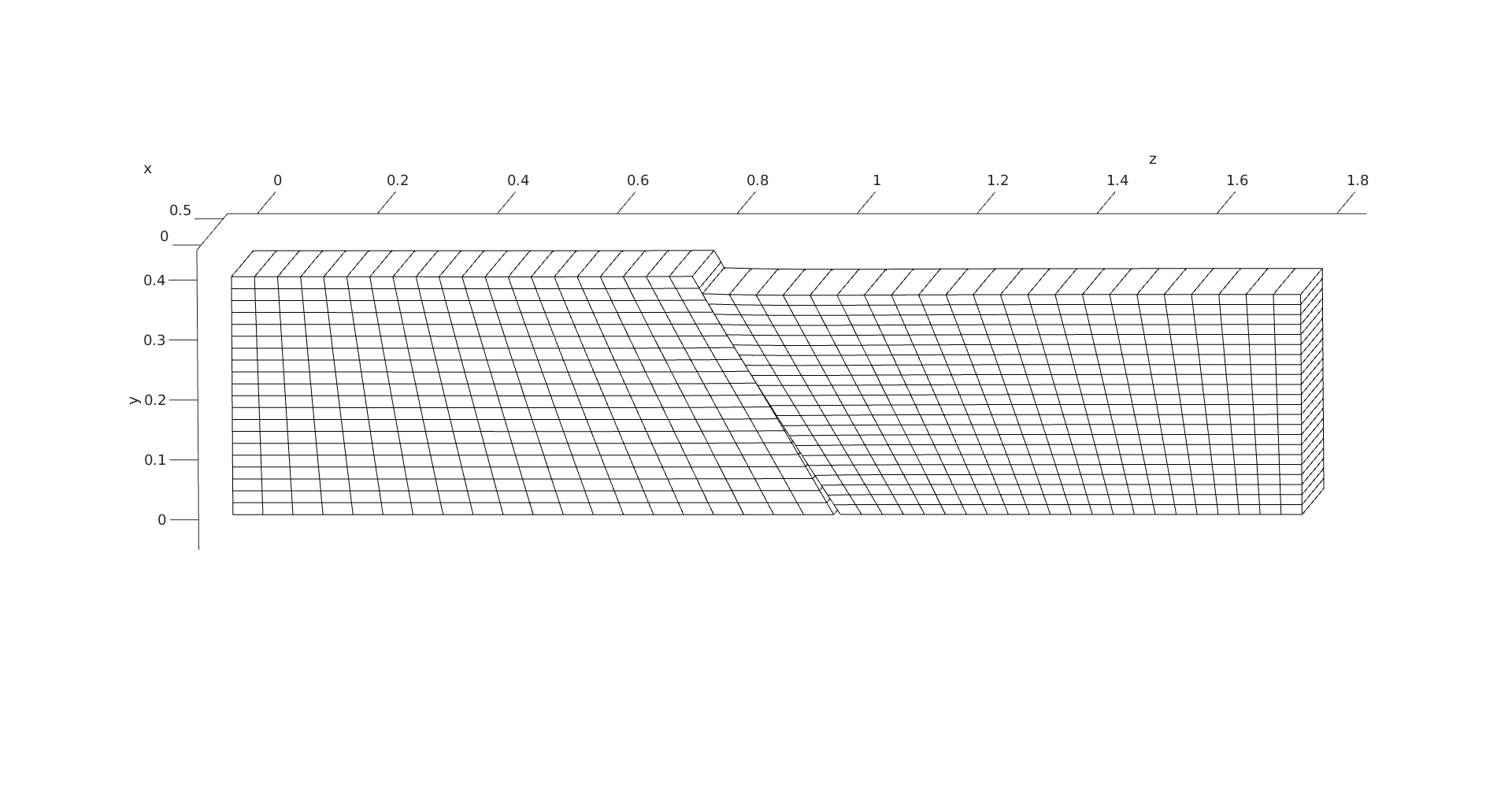}
     \end{subfigure} 
	\caption{The deformed configurations of the optimal design at time 5$s$ (top) and 10$s$ (bottom). } \label{fig:wedge-5s}
\end{figure}

\begin{figure}
  \centering
  \includegraphics[width=\textwidth]{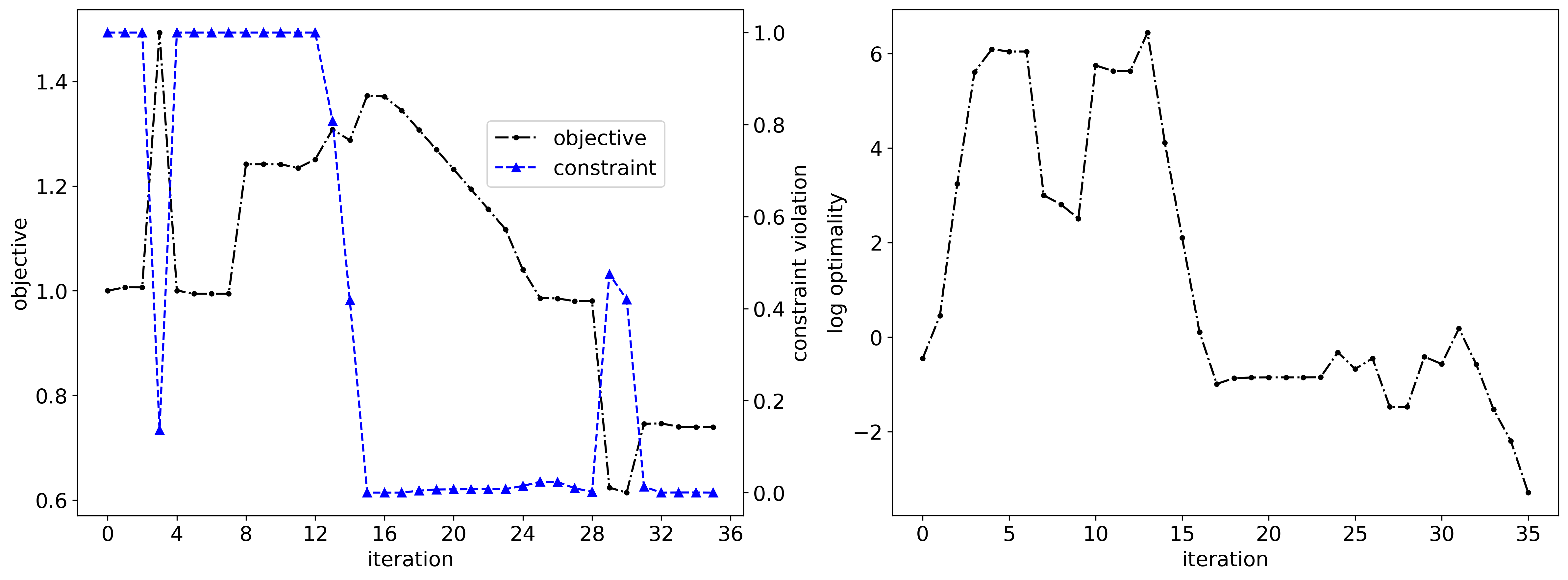}
	\caption{Optimization result for the wedge problem. The objective and constraint violation values are plotted on the left. The logarithm dual optimality measure is given on the right. }
\label{fig:wedge1}
\end{figure}

The initial values of the variables are set to $\theta_1=39^{\circ}$, $\theta_2=41^{\circ}$ and $P_1=1$.
Using the gradient-based approach, the optimal solution is successfully found in 35 iterations 
with the values $\theta_1 =  30, \theta_2 = 31.9583$ and $ P_1 =  0.7527$.
The deformed configuration of the optimal design at time $5$s and $10$s are shown in Figure~\ref{fig:wedge-5s}.
The optimization paths are plotted in Figure~\ref{fig:wedge1}.

Next, we report the results from constrained Bayesian optimization described in section~\ref{se:cbo}.  
We use Latin hypercube sampling to generate $8$ initial samples for the three dimensional design space. One more feasible sample point is added to start the algorithm.
The algorithm is given a total budget of $50$ iterations, comparable to the number of iterations required for convergence in the gradient-based approach. Thanks to the randomness of the algorithm, we repeat the same run $10$ times and present the average results.

\begin{figure}
  \centering
  \includegraphics[width=0.95\textwidth]{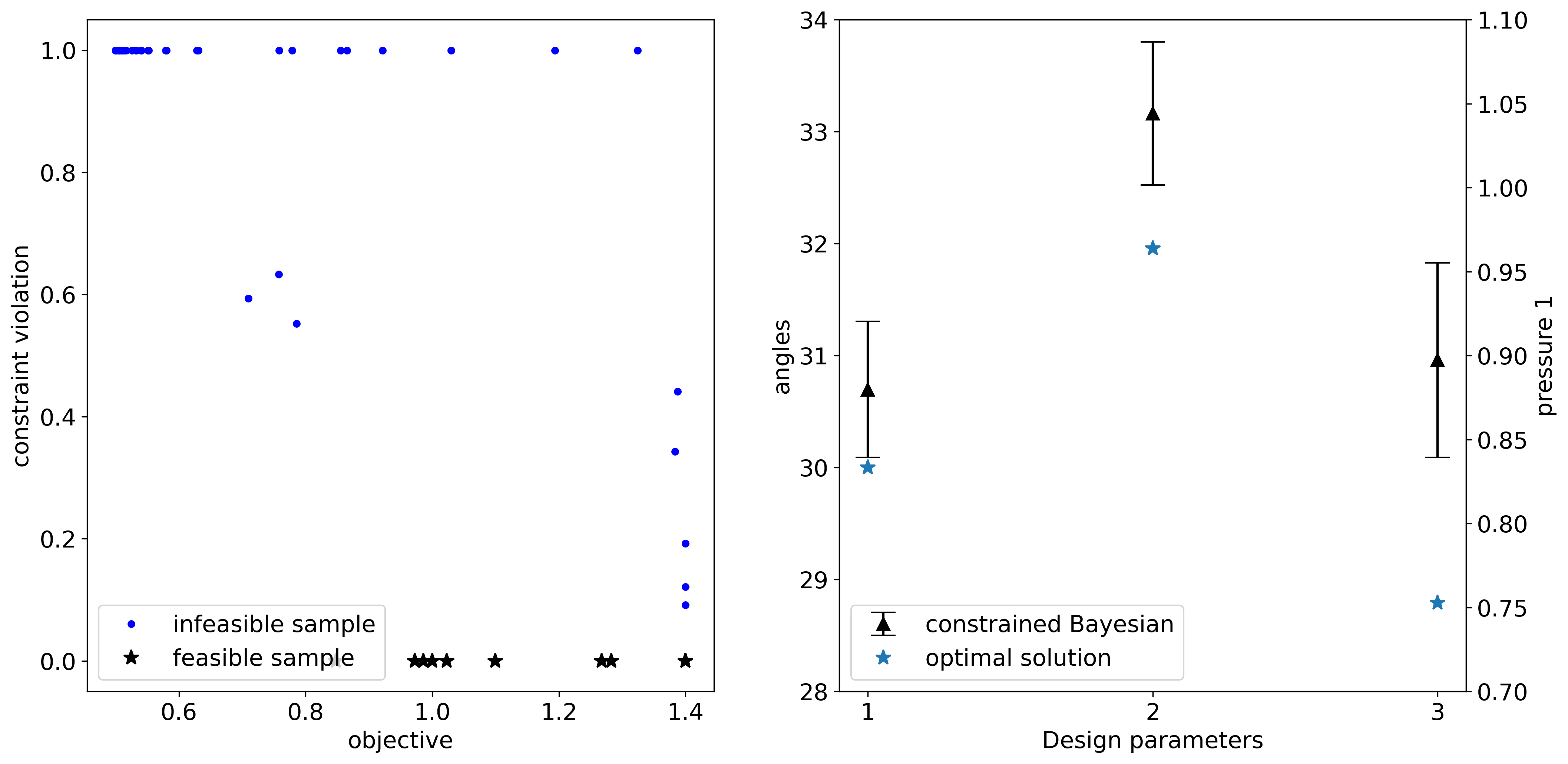}
	\caption{Constrained Bayesian optimization result for the wedge problem.}
\label{fig:cbo-wedge1}
\end{figure}

Since the feasible region is quite small, as evidenced by the high ratio of infeasible iterates during the optimization in Figure~\ref{fig:wedge1}, many samples are infeasible. However, the constrained Bayesian optimization succeeded in finding improved solution in $50$ iterations compared to the initial design. Figure~\ref{fig:cbo-wedge1} shows the objective and constraint violation of the samples collected during one of the runs, where the feasible and infeasible ones are marked with different signs. The optimal design variables are given on the right in Figure~\ref{fig:cbo-wedge1}, where the mean value and standard deviation of the ten constrained Bayesian optimization runs are computed and plotted, together with the optimal solution obtained through the gradient method. 
We note that throughout all the iterations and runs, the design obtained through gradient-based approach remains the optimal one.

\begin{figure}
  \centering
  \includegraphics[width=0.95\textwidth]{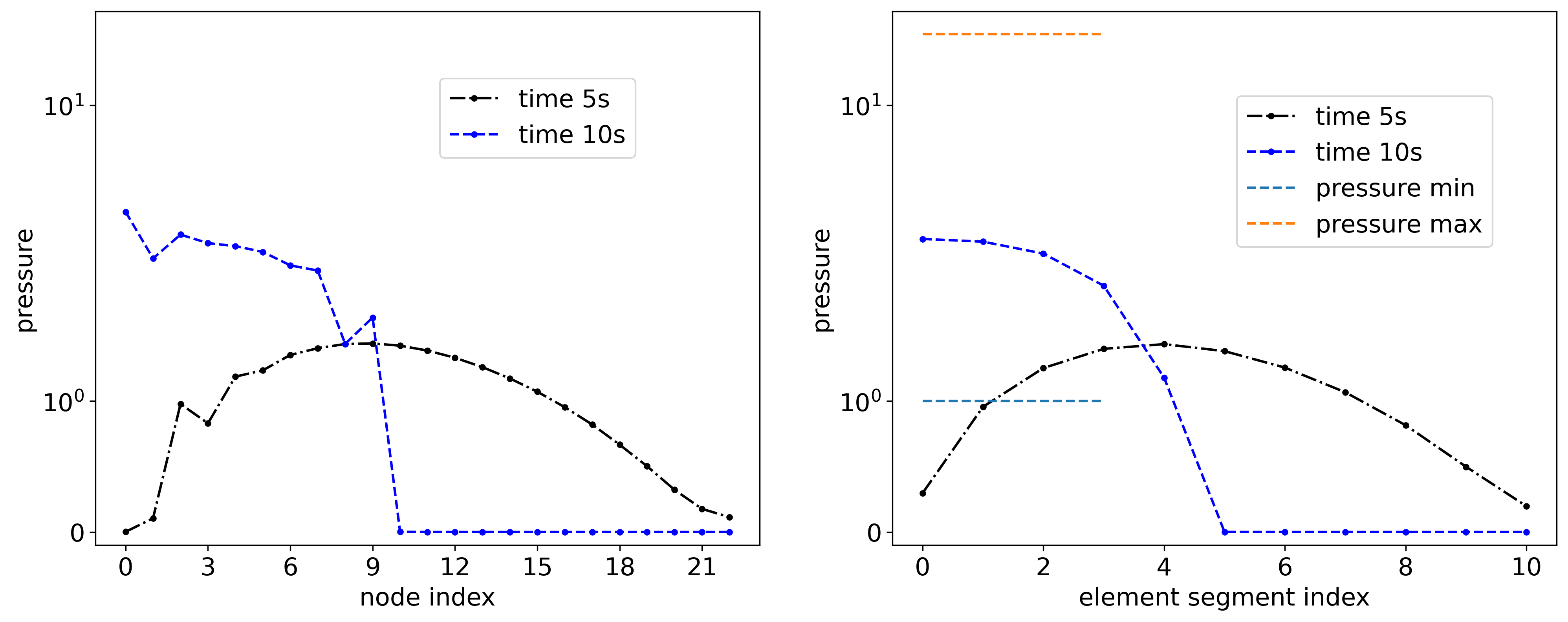}
	\caption{Pressure distribution on the contact surface for the initial design: nodal pressure (left) and averaged pressure on element segments (right). }
\label{fig:wedge-pressure1}
\end{figure}

\begin{figure}
  \centering
  \includegraphics[width=0.95\textwidth]{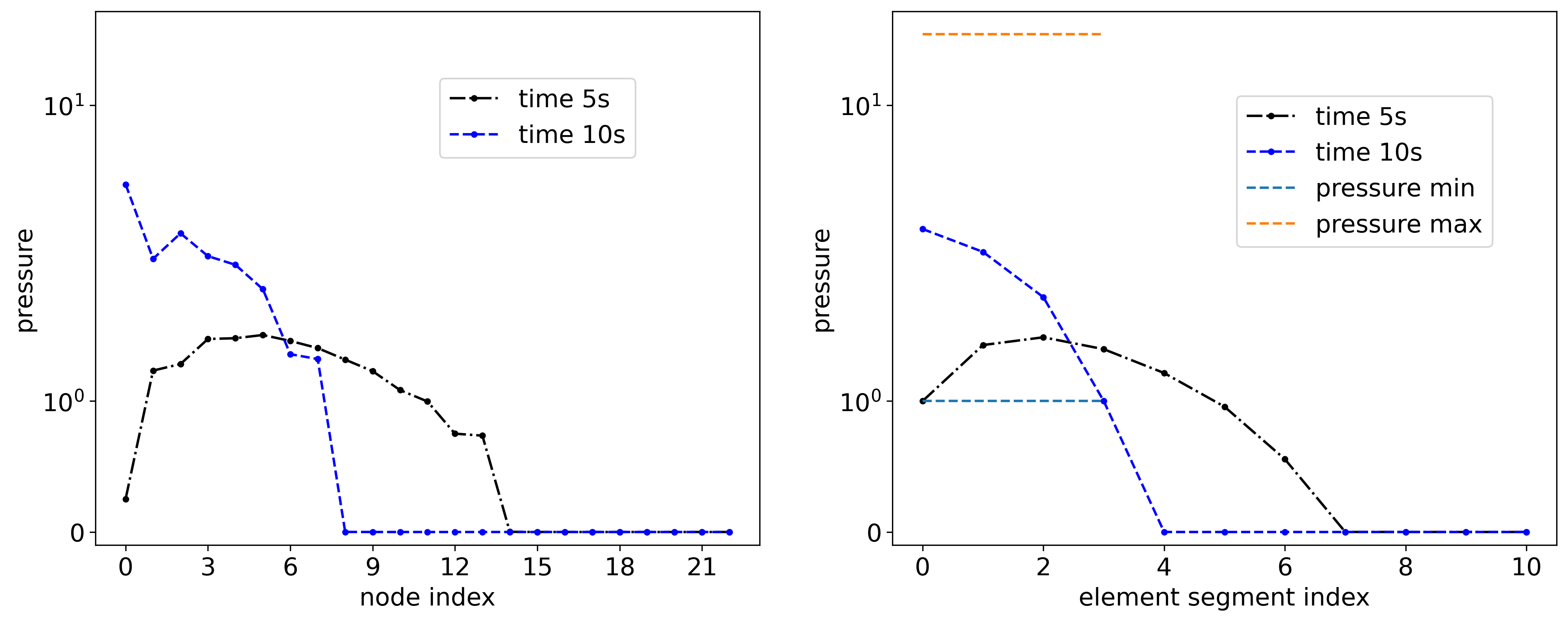}
	\caption{Pressure distribution on the contact surface for the optimal design: nodal pressure (left) and averaged pressure on element segments (right). }
\label{fig:wedge-pressure2}
\end{figure}

The optimization lays heavy emphasis on contact pressure, which is a critical source of nonsmoothness as discussed in previous sections. Hence, to better illustrate the feasibility region and how the pressure depends on optimization variables, we plot the pressure distribution of the $23$ nodes on the contact surface of the right wedge, as well as the averaged contact pressure of the $11$ element segments. The plot for the initial design $\theta_1=39,\theta_2=41,P_1=1$ is given in Figure~\ref{fig:wedge-pressure1}.	
The $x$-axis represents the nodes on the contact surface on the $x=0$ plane.
The order of increasing index matches with that of a decreasing $y$ coordinate, \textit{i.e.}, node $1$ is the upper-most node and element segment $1$ is the upper-most element segment on the contact surface.
The upper and lower bounds of pressure on the first four element segments are marked on the plots as well and the $y$-axis uses symlog scale for better illustration. 
The pressure plot for the optimal design $\theta_1 =  30, \theta_2 = 31.9583, P_1 = 0.7527$ is given in Figure~\ref{fig:wedge-pressure2}.
It is clear that the initial design does not produce a desired pressure distribution. 
Moreover, the averaging of pressure over element segments significantly smooths the nodal pressure distribution, which possibly leads to better convergence behavior of gradient-based optimization algorithms. 
Finally, we note that the two time steps of concern typically have quite different pressure distribution thanks to the different loading conditions. Therefore, satisfying pressure constraints at both time steps prove to be challenging.

\subsection{Marman clamp problem}
\begin{figure}
  \centering
  \includegraphics[width=0.85\textwidth]{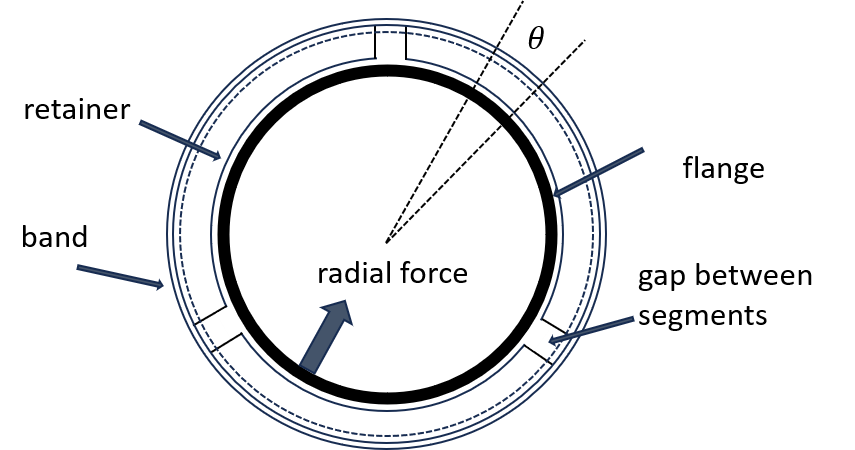}
	\caption{View looking down the cylindrical axis of a generic Marman clamp.  The flange is part of the cylindrical object that is being clamped together.  The gap between the segments of the retainer ensures that hoop stresses in the retainer are minimized.  The band is either stretched over the retainer such that the induced hoop stresses in the band induce compressive force on the retainer, or the compression may be induced by splitting the band at some point and pulling the pieces together using a screw mechanism. }
\label{fig:mclamp1}
\end{figure}

\begin{figure}
  \centering
  \includegraphics[width=0.44\textwidth,trim={0.1cm 0 0.2cm 0.cm},clip]{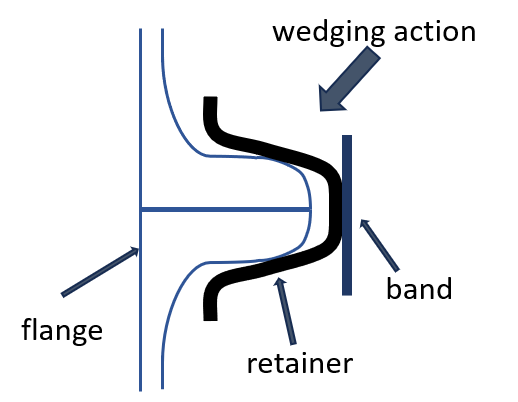}
  \includegraphics[width=0.55\textwidth,trim={0.cm 0 0.cm 0.cm},clip]{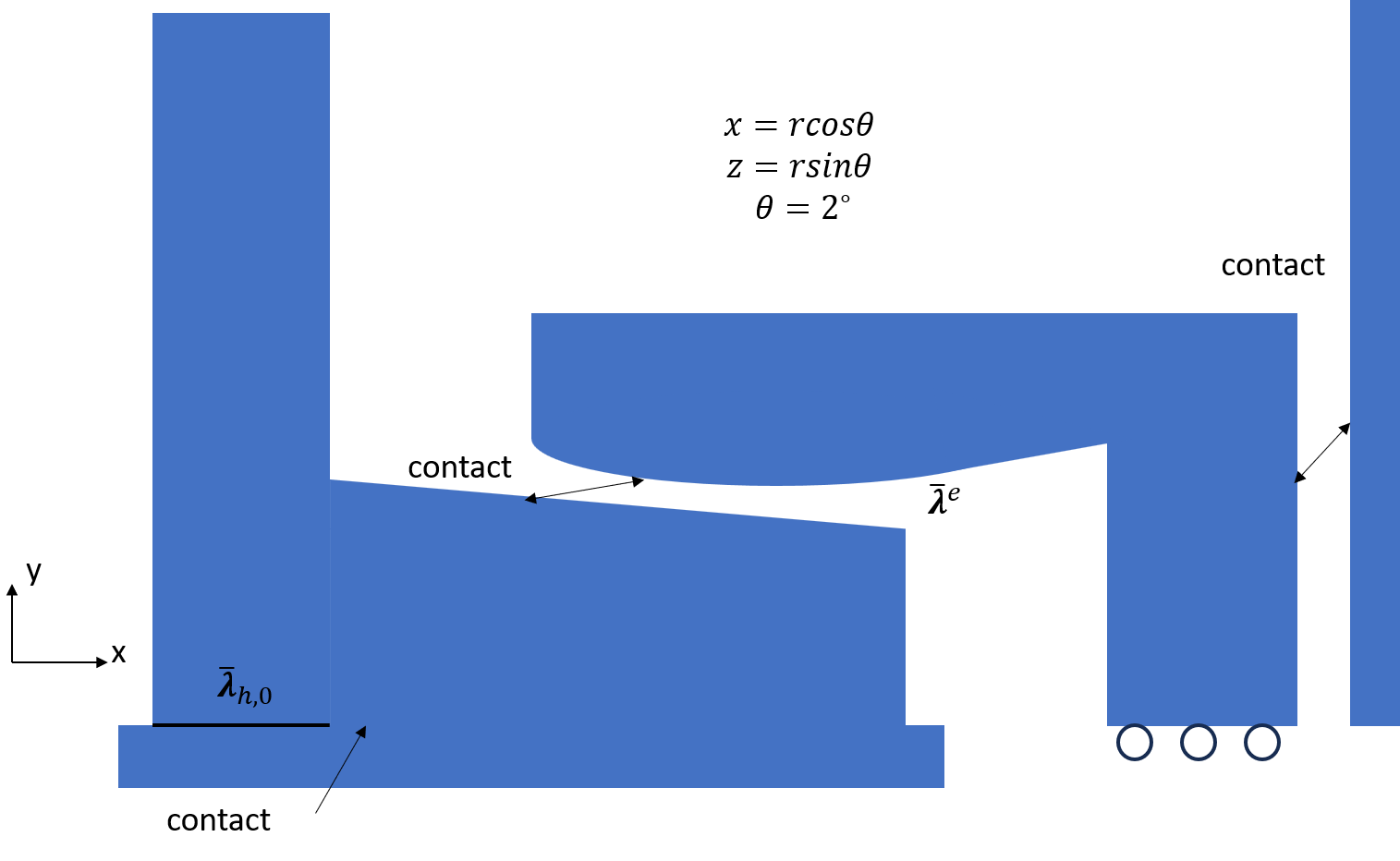}
	\caption{Left: Cross section of the clamp through a plane of constant $\theta$.  The cylinder is split into two pieces, the walls of the cylinder are terminated by flanges which meet at the centerline of the clamp. The retainer is forced against the flanges by the band, causing the flanges to be clamped against each other. Right: the contact surfaces and boundary conditions on an example design. The pressure involved in the constraints are marked.}
\label{fig:mclamp2}
\end{figure}

The second example considers the design of a Marman clamp~\cite{marmanclampguide}.  These devices are commonly used to connect two cylindrical sections together, see Figures~\ref{fig:mclamp1} and~\ref{fig:mclamp2}. In a general Marman clamp, the cylinders terminate in flanges which meet at a plane.  The two flanges interface with a retainer, which is forced against the flanges by the action of a clamp.  The retainers are commonly manufactured in segments to allow for inwards displacement without generating significant hoop stresses.  The shape of the retainer and the flanges forces the flanges together.  
In~\cite{yu2024lightweight}, the authors optimized clamp band joint using the so-called constrained sequential approximate optimization method, which bears some similarity to the constrained Bayesian optimization method. In~\cite{davar2020uncertainty}, uncertainty analysis for the clamp band joints is incorporated into a mathematical model.

In the manifestation considered here, the clamp is in the form of a band which is forced against the retainer by designing the components such that there is an interference fit between them.  The hoop stress generated in the band is reacted by a contact pressure between the band and the retainer, which is further reacted by contact pressure between the retainer and the flanges, and ultimately the two flanges are forced together, inducing contact pressure between them.

The cylindrical structures can be treated as axisymmetric so that only a small portion of the clamp is simulated and the angular coordinates ($\theta$) are fixed in the cylindrical coordinate system. Given a $\theta$, the bodies in the radial and longitudinal plane, $x-y$ plane, is shown in Figure~\ref{fig:mclamp2}.
Similarly, only the top flange (referring to Figure~\ref{fig:mclamp2}) and the top half of the band and the clamp are simulated.  Three contact planes are defined – contact between the band and the retainer, between the retainer and the top flange, and between the top flange and the vertical symmetry plane, represented by an artificial box which is held rigid.  The vertical displacement at the vertical symmetry plane of the retainer and the band is held fixed in accordance with typical symmetry boundary conditions.
As mentioned earlier, hoop stress in the band is generated through an initial interference fit, represented as an initial overlap in the meshes.

Both the flange and the band are given isotropic linear elastic materials with the elastic modulus at $E = 2\times 10^6$ and Poisson's ratio at $\nu =0.3$. 
The retainer is given an anisotropic material that has a smaller ($0.1\%$) elastic and shear modulus in the $\theta$-direction in order to simulate the gap between segments of the retainer (Figure~\ref{fig:mclamp1}).

\begin{figure}
  \includegraphics[width=0.5\textwidth,trim={28cm 2cm 28cm 0cm},clip]{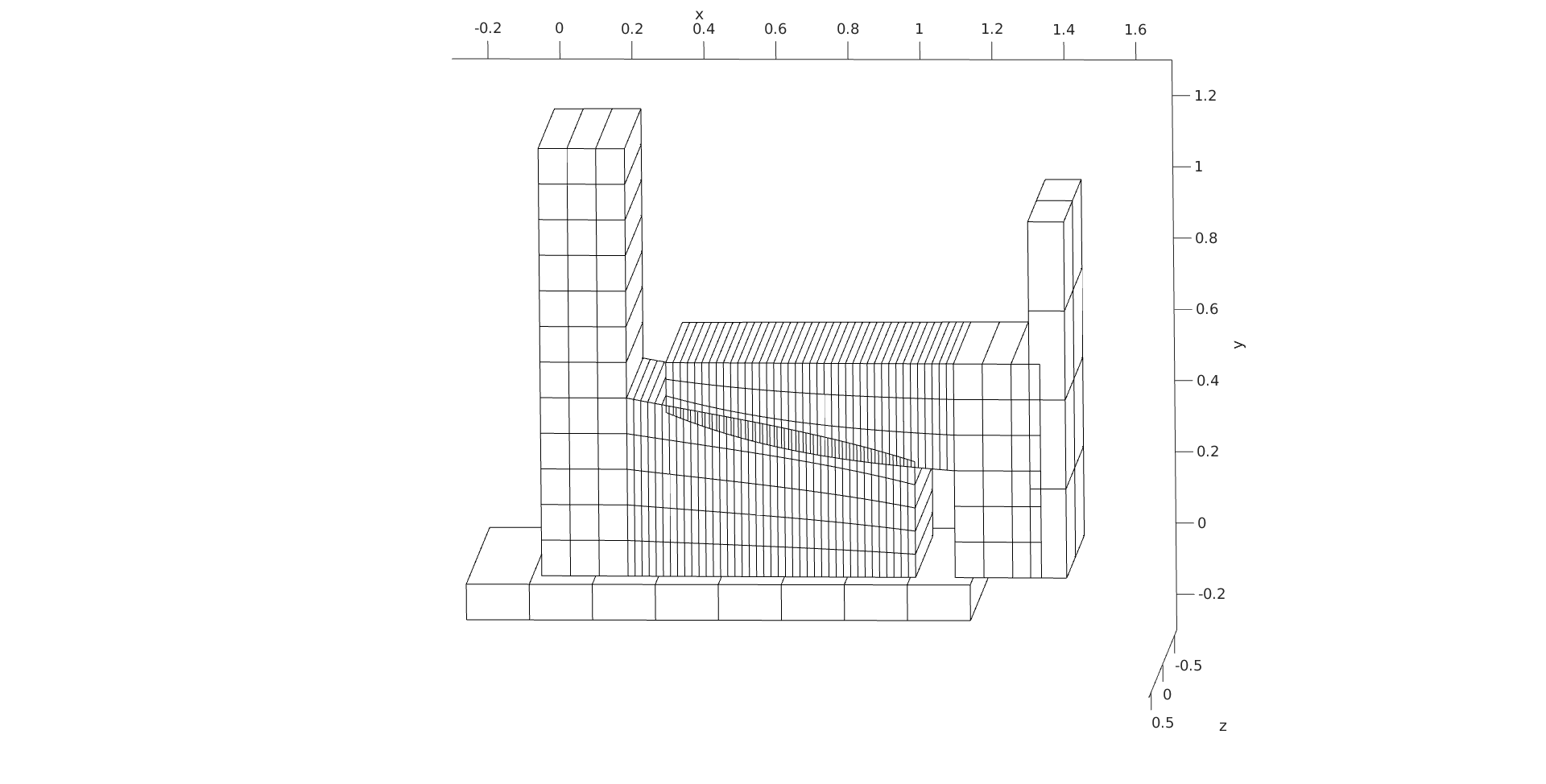}
  \includegraphics[width=0.5\textwidth,trim={28cm 2cm 28cm 0cm},clip]{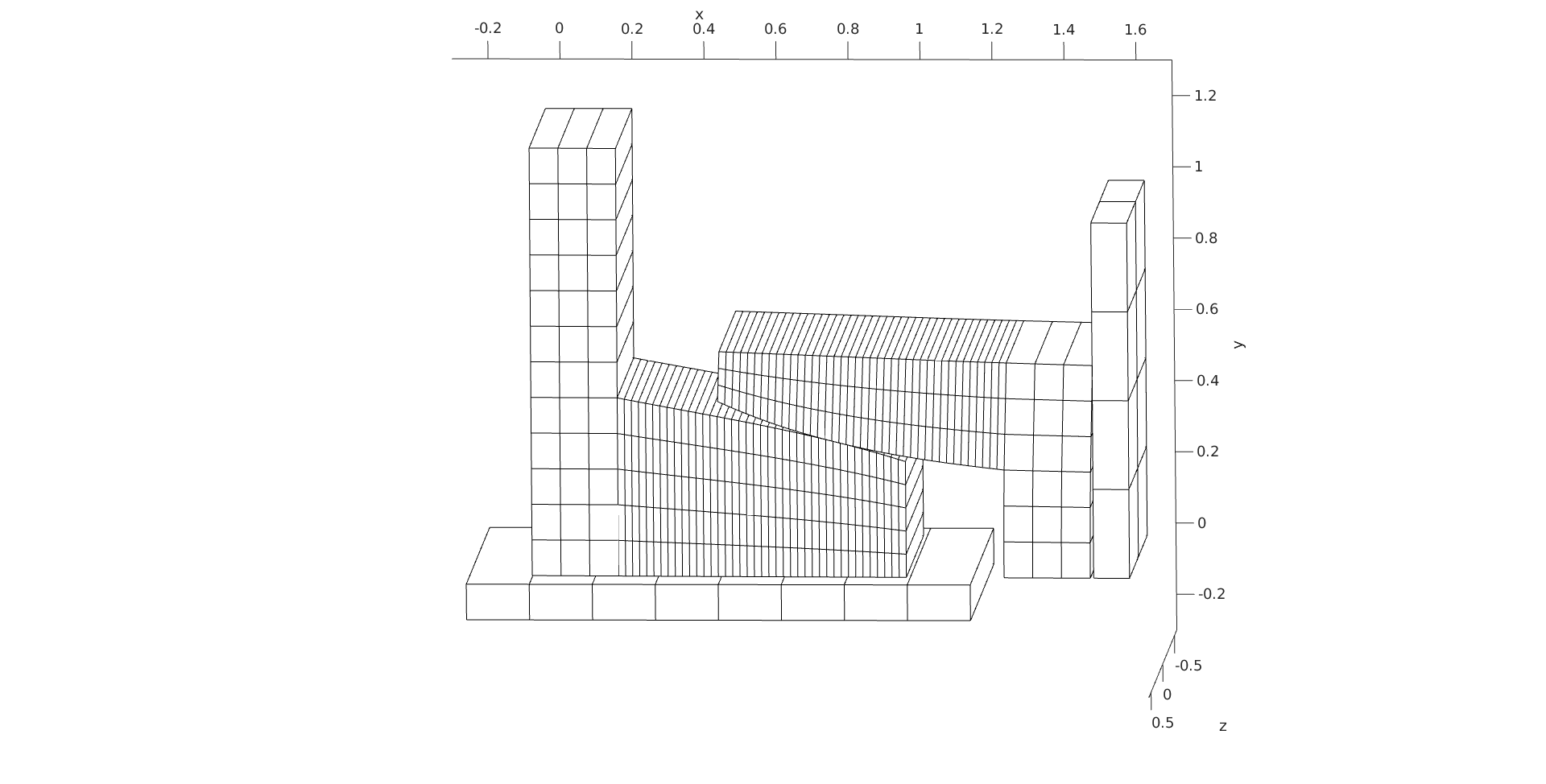}
	\caption{Finite element mesh of the initial design in the undeformed configuration (left) and deformed configuration (right).  Note that the band is meshed with an initial penetration which will be forced out by the contact, inducing hoop stresses in the band and compressive stresses in the retainer.  The bodies are meshed as a 3D wedge with angled boundary conditions to enforce axisymmetry.}
\label{fig:mclamp-init}
\end{figure}

The referential three-dimensional finite element model of an example design (the initial design) is shown on the left in Figure~\ref{fig:mclamp-init}. 
The mesh around the potential contact surface between the retainer and the flange is refined so that a more accurate contact solution can be obtained.
Notice that the initial configurations have penetrations, which are resolved as part of the contact forward solution.

To ensure robustness of the design against separation of the flange, the objective is chosen to be minimization of the compliance, that is maximizing the stiffness of the overall system.
Meanwhile, the maximum contact pressure between the retainer and the flange is constrained by user-defined upper and lower bounds. 
The maximum pressure is imposed on the element pressure to reduce the nonsmoothness of the nodal pressure distribution. 
The element pressure is denoted as $\bar{\lambdabold}^e$ and computed as the average of its nodal pressure . 
The $40$ elements on the contact surface are used to define the constraints. 
To prevent leaks from the cylindrical assembly and ensure adequate sealing, we enforce an interface force at the bottom-left part of the flange.
This is reflected as a constraint on the sum of the four nodal pressure at the bottom-left part of the flange, given that the mesh and surface area of this region is independent of the design variable.

\begin{figure}
  \centering
  \includegraphics[width=0.85\textwidth]{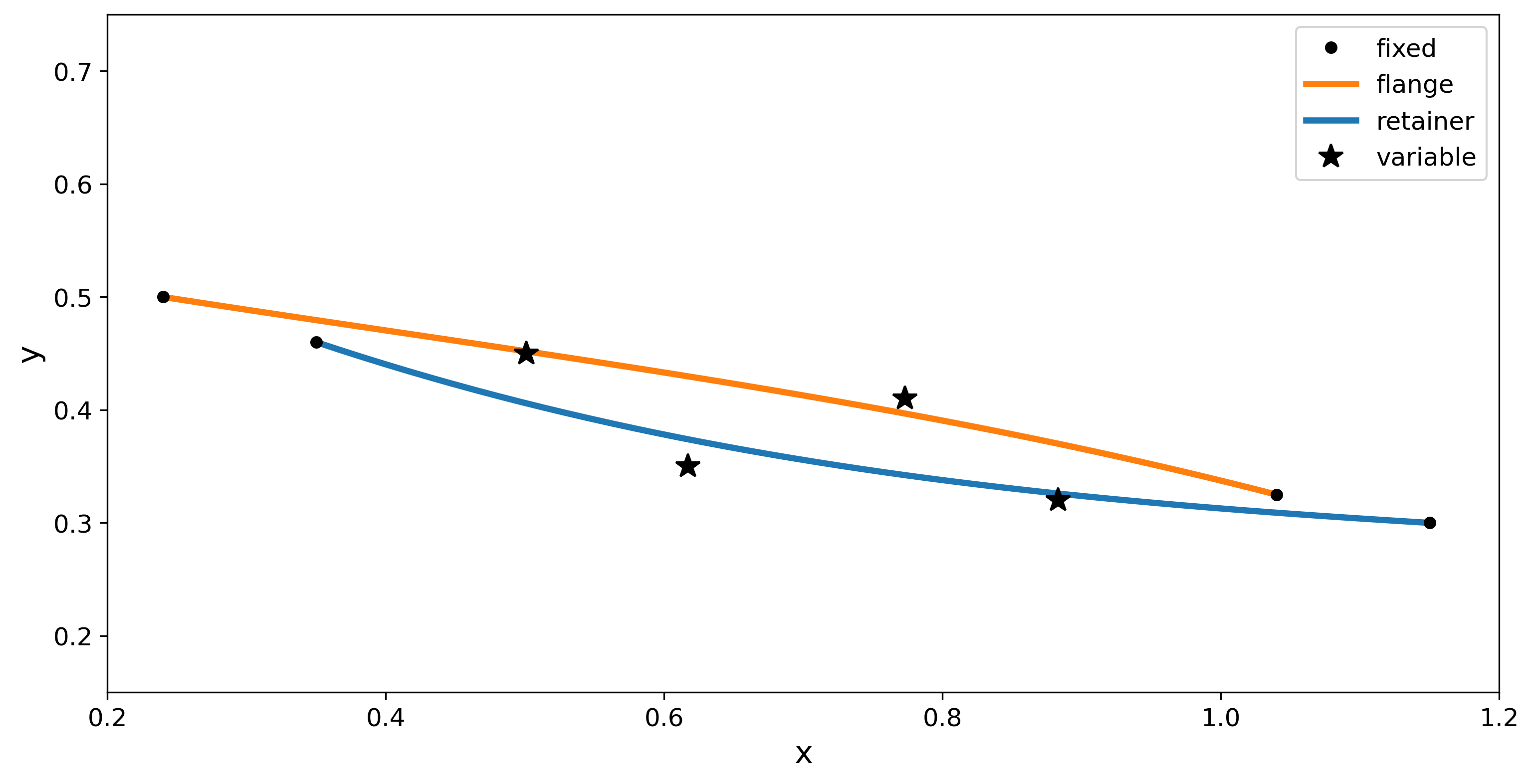}
	\caption{The curve to be optimized in the $x-y$ plane. The $y$-coordinates of the `variable' nodes are design variables. }
\label{fig:bezier_curve_2d}
\end{figure}

The design variables are the contact geometries of the retainer and flange. Specifically, in the $x-y$ plane the contact surfaces are defined by two cubic Bézier curves, on both the retainer and the flange, an example of which is shown in Figure~\ref{fig:bezier_curve_2d}. Given four nodes $\Pbm_i,i=1,2,3,4$, the explicit form of a curve is
\begin{equation} \label{def:cubic-bezier}
 \centering
  \begin{aligned}
    \Bbm(t) = (1-t)^3 \Pbm_1+3(1-t)^2 t \Pbm_2 + 3(1-t)t^2 \Pbm_3+t^3\Pbm_4, \ 0\leq t\leq 1, 
  \end{aligned}
\end{equation}
where the vector $\Bbm\in\Rbb^3$ is the nodal coordinate.
The curve passes through $\Pbm_1$ and $\Pbm_4$ but not necessarily the two middle control points $\Pbm_2$ and $\Pbm_3$.

Given the geometry and assembly constraints, for both the retainer and the flange, the two end points $\Pbm_1$ and $\Pbm_4$ of the two curves are fixed, as well as the $x$-coordinates of the two middle control points $\Pbm_2$ and $\Pbm_3$ of each curve.
The $y$-coordinates of the four middle control points, two from each curve, are the optimization variables. Their $z$-coordinates are determined by cylindrical symmetry, where $\theta=2^{\circ}$ is fixed.
The inner radius of the cylinder as seen in Figure~\ref{fig:mclamp1} is $20$.
Further, each variable is bounded in order to ensure valid retainer and flange shapes as well as contact between the retainer and flange. The bounds are simply chosen through numerical experiments and observation.  

The mathematical form of the optimization problem is 
\begin{equation} \label{eqn:opt-1st-clamp}
 \centering
  \begin{aligned}
   &\underset{\substack{\rhobold \in \Rbb^4}}{\text{minimize}} 
	  & &  a(\rhobold, \bar{\ubm}_h)\\
   &\text{subject to}
	  & & \lambda_{l,o}\leq \bar{\lambda}_{h,o} \\
          &&& \bar{\lambda}_l \leq \max (\bar{\lambdabold}^e) \leq \bar{\lambda}_{u}\\
          &&&   \rhobold_l \leq \rhobold \leq \rhobold_u,\\
  \end{aligned}
\end{equation}
where $\rhobold \in\Rbb^4$, $\bar{\lambdabold}^e$ is the pressure vector of the element segments on the contact surface between the flange and the retainer. The objective function $a$ is the compliance function. The user-defined constraint bounds heavily influence the optimal design and can be tuned for other design objectives.
The pressure  $\bar{\lambda}_{h,o}$ is the sum of pressure of the four nodes on the bottom-left part of the flange.
In this example, we set the $x$-coordinates of the four middle control points to $[0.617,  0.883]$ for the retainer and $[0.50,  0.773 ]$ for the flange so that all four control points are evenly distributed in the $x$-direction.
We choose the bounds to be $\rhobold_u =[0.408, 0.354, 0.47, 0.44]$ and $\rhobold_l = [0.35,0.35,0.443, 0.3834]$ and constraint parameters to be $\lambda_{l,o}=30$, $\bar{\lambda}_l = 300$ and $\bar{\lambda}_u = 650$.

The initial $\rhobold$ is set to $[0.35, 0.32,  0.45, 0.41]$. Its deformed configuration from contact simulation is shown on the right in Figure~\ref{fig:mclamp-init}.
The pressure on the bottom-left nodes of the flange is essentially $0$, as the contact between the flange and the retainer occurs towards the right end of the flange. 
For the gradient-based approach, an optimal solution is successfully found in 48 iterations with Ipopt, using the adaptive barrier update strategy to a dual optimality accuracy of $10^{-4}$. 
The optimal set of values are $[0.40713,    0.35344,  0.443, 0.38343]$, with three of the four variables not reaching bounds.

\begin{figure}
  \centering
  \includegraphics[width=1.0\textwidth]{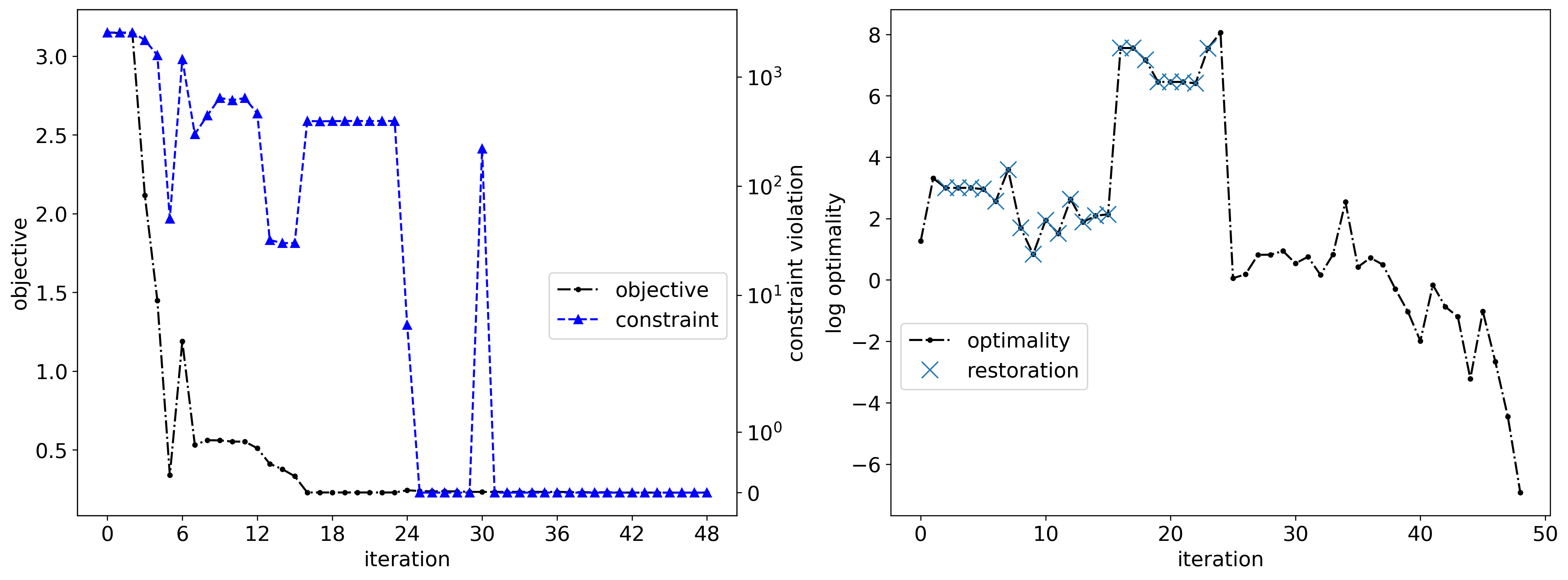}
	\caption{The optimization iterations for the clamp problem using Ipopt. The objective and constraint violation values are given on the left. The logarithm of dual optimality measure is given on the right. }
\label{fig:mclamp-opt}
\end{figure}

\begin{figure}
  \centering
  \includegraphics[width=0.49\textwidth,trim={35cm 2cm 30cm 0cm},clip]{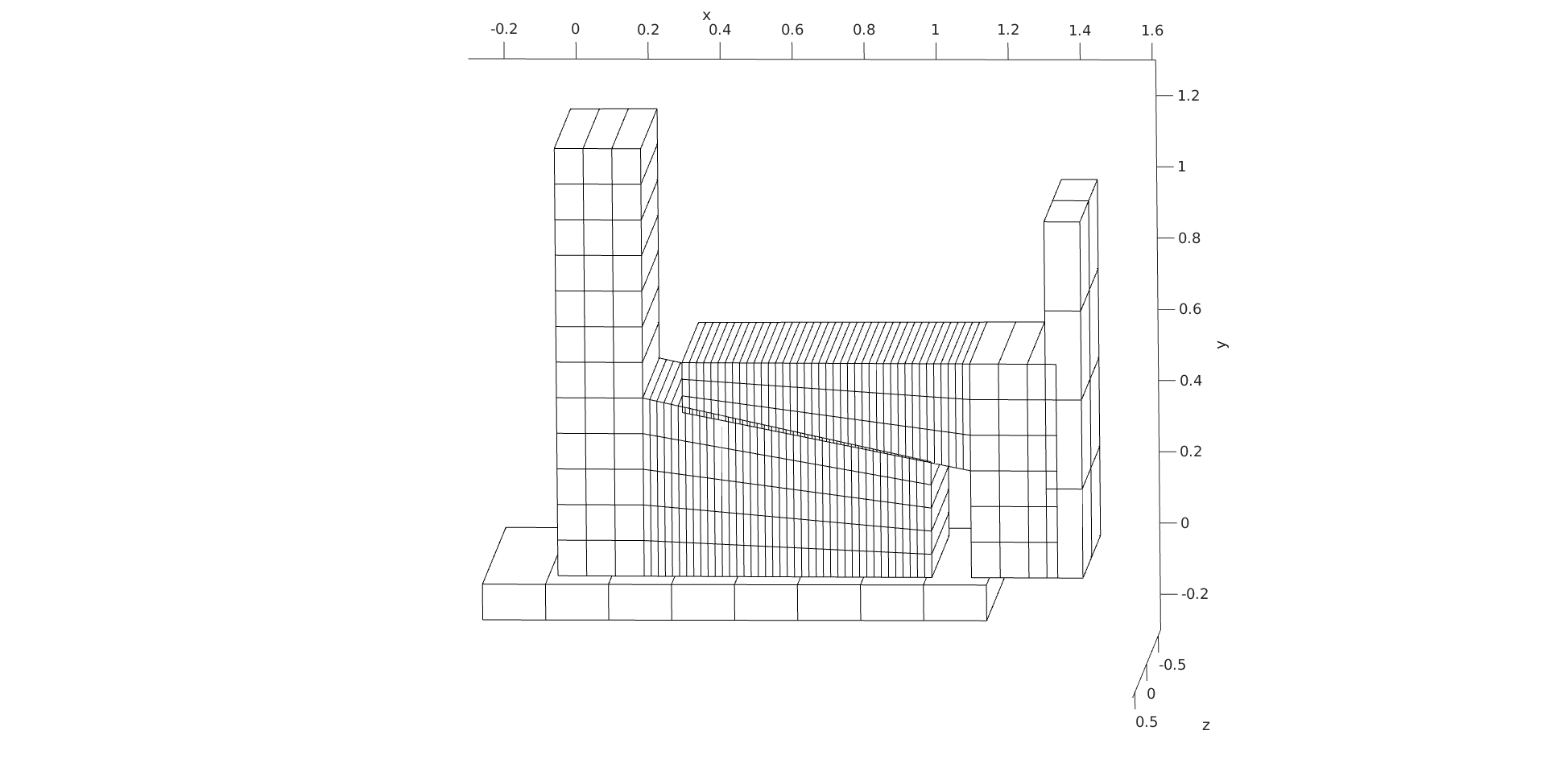}
  \includegraphics[width=0.49\textwidth,trim={35cm 2cm 30cm 0cm},clip]{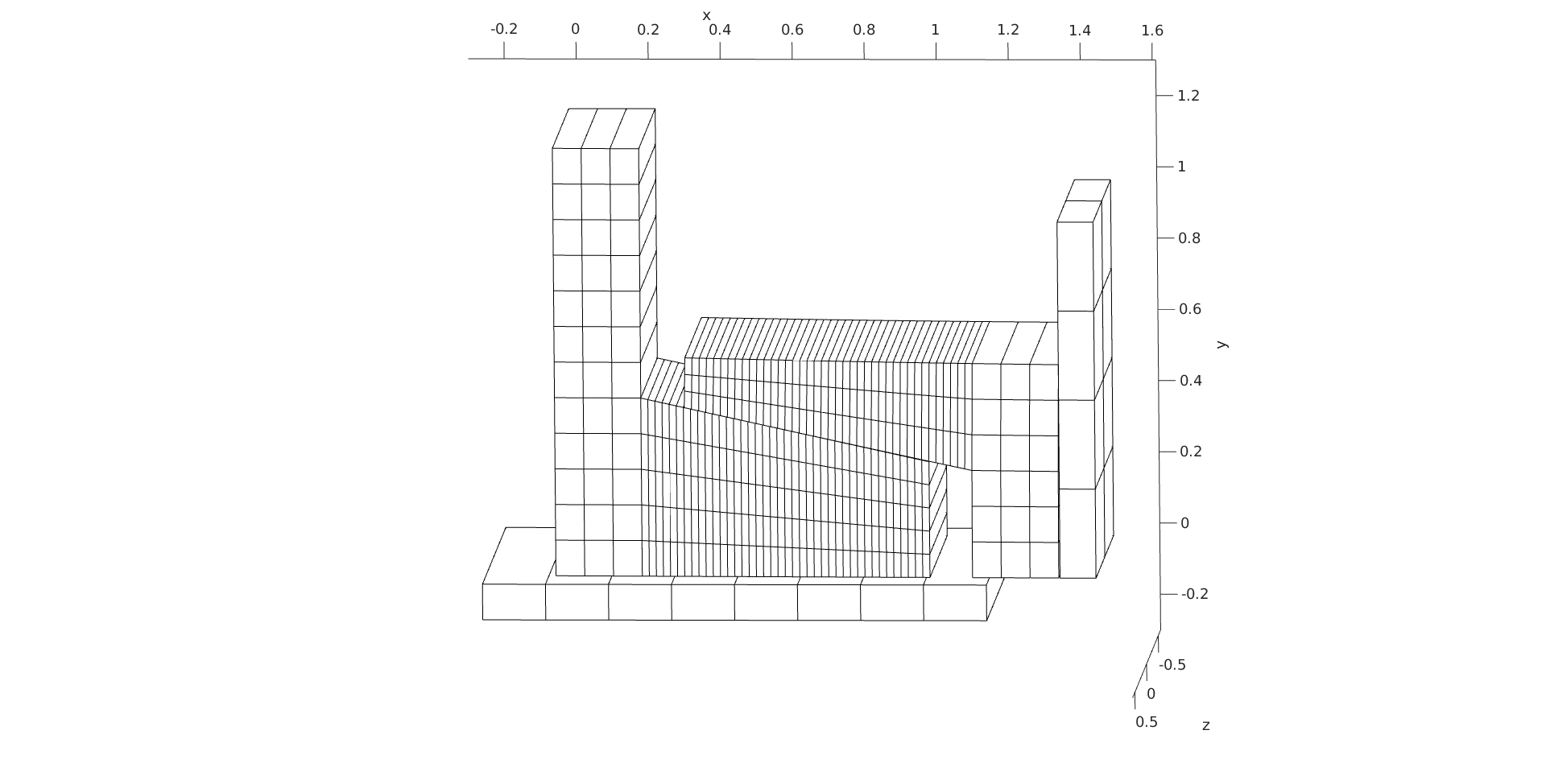}
	\caption{The undeformed (left) and deformed (right) configurations of the optimal design. }
\label{fig:mclamp-optimal}
\end{figure}

\begin{figure}
  \centering
  \includegraphics[width=1.0\textwidth]{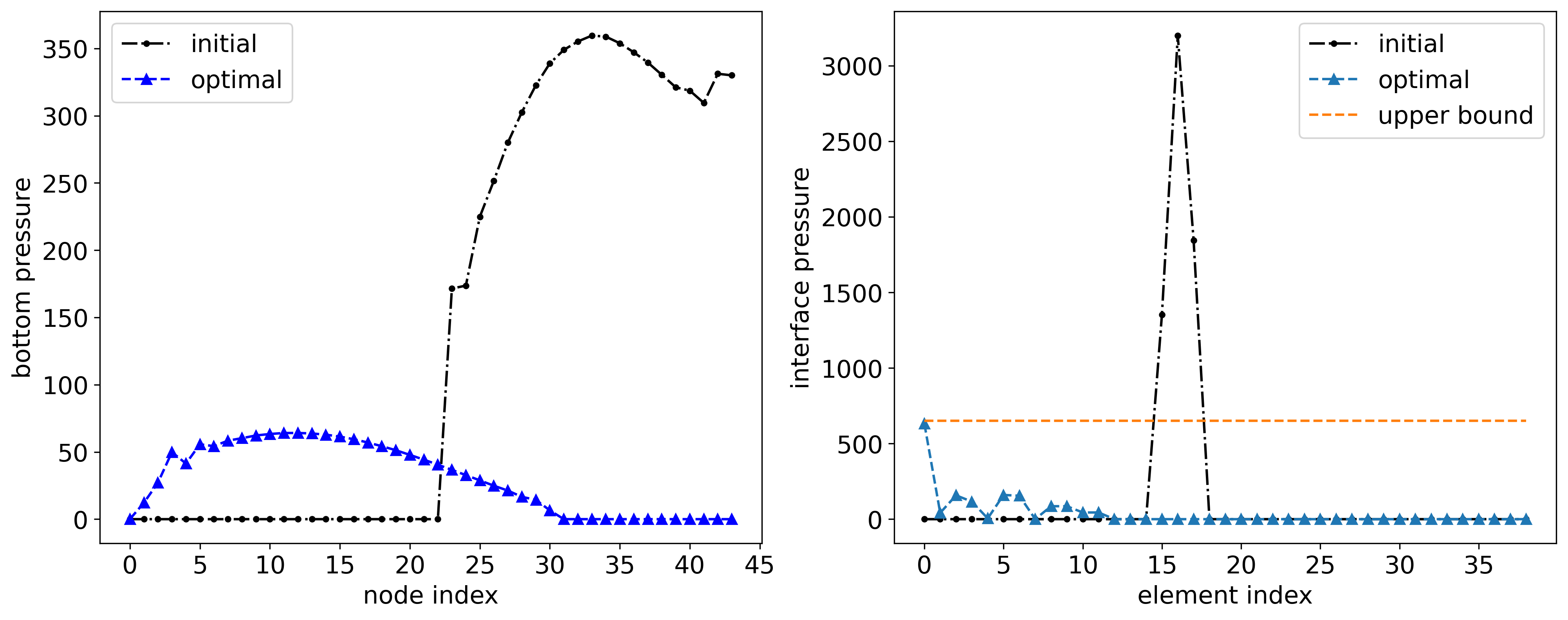}
	\caption{Pressure distribution of the initial and optimal design at the bottom contact surface (left) and the contact surface between the flange and the retainer. The node index for the bottom contact surface increases with the $x$-coordinate. The element index of the contact surface between the flange and the retainer increases with the $x$-coordinate as well.}
\label{fig:mclamp-pressure}
\end{figure}

The objective, constraint violation and dual optimality during the optimization are plotted in Figure~\ref{fig:mclamp-opt}, where the $y$-axis of constraint violation uses symlog scale for better illustration. Additionally, the feasibility restoration steps are marked with the $\times$ sign. 
The deformed configuration of the optimal design is shown in Figure~\ref{fig:mclamp-optimal}. 
Compared to the initial design, the contact region has moved to the left so that pressure could be passed onto the bottom-left part of the flange. 
The pressure distribution of the bottom contact surface and the contact surface between the flange and the retainer is plotted in Figure~\ref{fig:mclamp-pressure}. For the contact surface between the flange and the retainer, we plot the element pressure. It is worth noting that thanks to the curved design of the two contacting surfaces, the element pressure on the contact surface between the flange and the retainer is not smooth or monotonic. 

While the contact nodes and regions change with the design variables, sensitivities with Ipopt manage to overcome potential nonsmoothness to find a reasonable solution through careful design of the bound constraints. Future work to better understand the nonsmoothness and design more robust optimization algorithms is being planned. 

Next, we apply the constrained Bayesian optimization algorithm (Algorithm~\ref{alg:bo}) to the clamp design problem, where $16$ initial Latin hypercube samples are used due to the increased dimension of the design space. One more feasible sample point is added to start the algorithm. We run the Bayesian algorithm for $50$ iterations and repeat the same run $10$ times to obtain less uncertain results.  

As evidenced by the large number of feasibility restoration steps with Ipopt, the feasible region for this problem is rather small.
Applying BFGS interior-point method to find the maximum of the constrained expected improvement acquisition function (line 4 of Algorithm~\ref{alg:bo}) is not effective due to the low probability of constraints being satisfied.
Therefore, we use random search method for this problem where the maximum constrained expected improvement value is chosen from $10000$ random variables at each iteration to determine the next sample. We note that only the constrained expected improvement function, which has a closed form, is evaluated at these variables. No additional contact forward simulation is required. 

\begin{figure}
  \centering
  \includegraphics[width=0.6\textwidth]{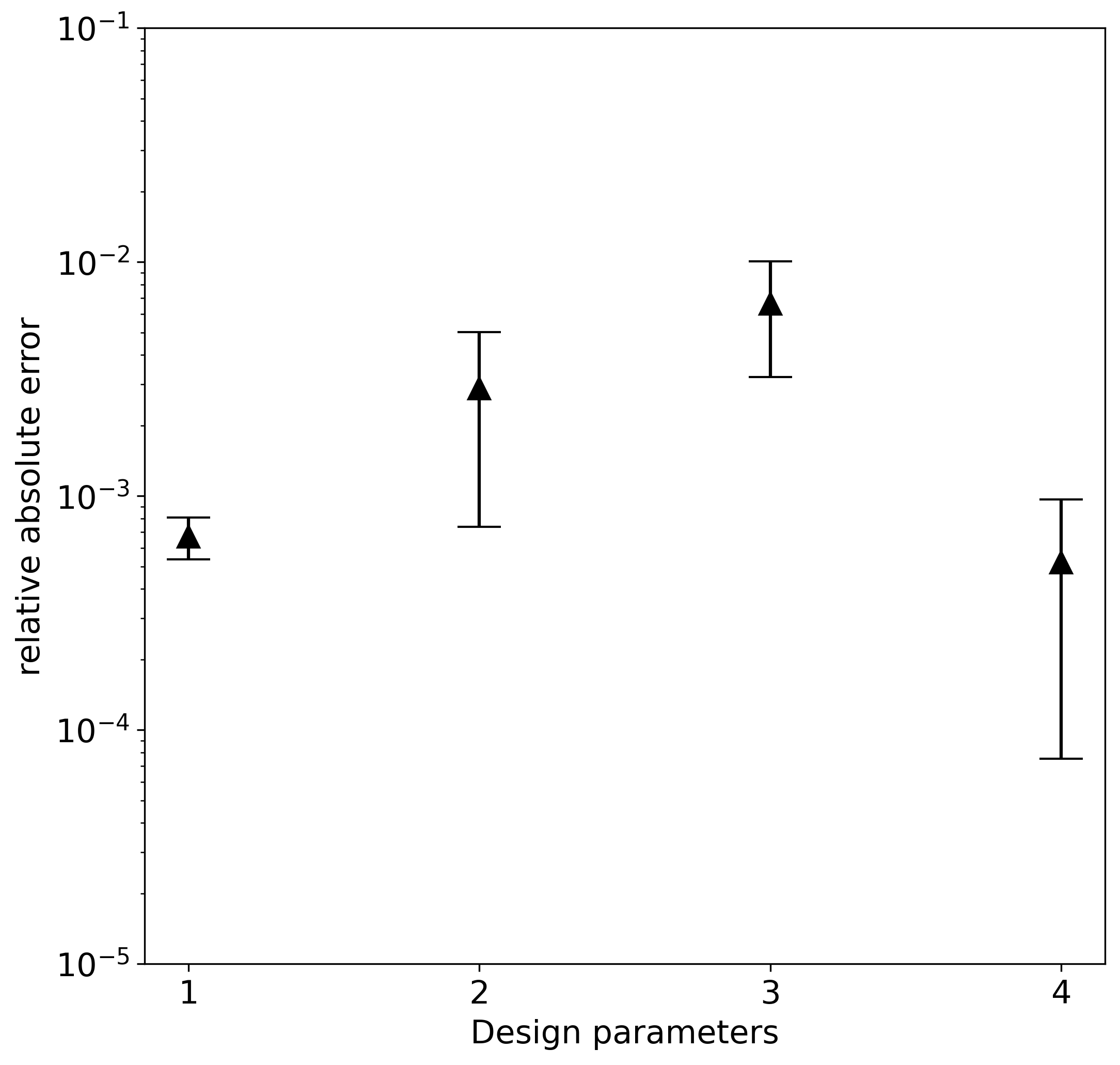}
	\caption{The relative absolute error of the mean optimal variables from constrained Bayesian optimization. Plot is in log scale.}
\label{fig:mclamp-baye}
\end{figure}

On average, we obtain $6$ out of the $66$ total samples as feasible designs. We compute the relative absolute error of the mean value of the optimal design from Bayesian optimization from the optimal design from Ipopt, as plotted in~\ref{fig:mclamp-baye}. Overall, the Bayesian optimization algorithm generates designs close to the one from the gradient-based method. 

The overall results from the constrained Bayesian optimization show that an improved design can be obtained for a reasonable sample size. Compared to the gradient-based approach, it could offer an alternative way for optimization with higher efficiency but reduced accuracy, particularly for design problems with a relatively small number of variables and no available derivatives.  
However, the constrained Bayesian optimization algorithm adopted here requires a feasible initial sample, which could limit its applicability. The authors are proposing a new and improved constrained Bayesian optimization method in~\cite{wang2024constrained}. 
We also point out that for the general topology optimization problems where the densities are the variables, the large dimension of the problem could make Bayesian optimization impractical.

\section{Conclusion}\label{se:conclusion}
In this paper, we solved design optimization problems in unilateral contact under pressure constraints, using both gradient-based optimization approach and a gradient-free constrained Bayesian optimization algorithm. 
The contact finite element simulation uses the mortar gap function formulation and interior-point forward solver.
Through smoothing of nodal pressure and careful design of the bound constraints, we successfully find optimal solutions using the gradient-based method and sensitivities, even though the problems are considered nonsmooth. 
The Bayesian optimization approach produced relatively good results and could be suitable for problems without reliable sensitivity information. To solve our target problems more efficiently and robustly, we should address the nonsmoothness directly and develop algorithms that might overcome it. We aim to work on this challenge in the future.

\section*{Acknowledgments}
The authors would like to thank Professor Daniel A. Tortorelli for his help and support in preparation of this manuscript.
The valuable feedback from Professor Tortorelli has contributed to the manuscript. 
Prepared by LLNL under Contract DE-AC52-07NA27344. Release number: LLNL-JRNL-863659-DRAFT. 
\bibliography{bibliography}

\end{document}